\documentclass[aap,MSNbibl,citesort,dvips]{arximspdf}
\usepackage{accents}
\usepackage{graphicx}

\doi{10.1214/11-AAP777}
\volume{22}
\issue{1}
\pubyear{2012}
\firstpage{407}
\lastpage{455}

\makeatletter

\renewcommand{\underbar}{\underaccent{\bar}}

\newcommand{\Pd}{\mathbb{P}}
\newcommand{\akl}{a^{k_l}}

\newcommand{\bbR}{\mathbb{R}}

\newcommand{\bbZ}{\mathbb{Z}}
\newcommand{\Ex}{\mathbb{E}}

\newcommand{\mK}{\mathcal{K}}
\newcommand{\mH}{\mathcal{H}}
\newcommand{\mB}{\mathcal{B}}
\newcommand{\mC}{\mathcal{C}}

\newcommand{\argmax}{\arg\max}
\newcommand{\argmin}{\arg\min}

\newcommand{\lam}{^n}

\newcommand{\mN}{\mathcal{N}}

\newcommand{\mF}{\mathcal{F}}
\newcommand{\hX}{\hat{X}}
\newcommand{\hW}{\hat{W}}

\newtheorem{theorem}{Theorem}[section]
\newtheorem{prop}{Proposition}[section]
\newproclaim{defin}{Definition}[section]

\newproclaim{rem}{Remark}[section]
\newtheorem{lem}{Lemma}[section]

\newcommand{\I}{\mathcal{I}}
\newcommand{\tinf}{\rightarrow\infty}

\makeatother

\begin{document}
\begin{frontmatter}

\title{On optimality gaps in the Halfin--Whitt regime}
\runtitle{Optimality gaps in the Halfin--Whitt regime}

\begin{aug}
\author[A]{\fnms{Baris} \snm{Ata}\ead[label=e1]{b-ata@kellogg.northwestern.edu}} and
\author[A]{\fnms{Itai} \snm{Gurvich}\corref{}\ead[label=e2]{i-gurvich@kellogg.northwestern.edu}}
\runauthor{B. Ata and I. Gurvich}
\affiliation{Northwestern University}
\address[A]{Kellogg School of Management\\
Northwestern University\\
Evanston, Illinois 60208\\
USA\\
\printead{e1}\\
\hphantom{E-mail: }\printead*{e2}}
\end{aug}

\received{\smonth{6} \syear{2010}}
\revised{\smonth{12} \syear{2010}}

%
\begin{abstract}
We consider optimal control of a multi-class queue in the Halfin--Whitt
regime, and revisit the notion of asymptotic optimality and the
associated optimality gaps. The existing results in the literature for
such systems provide asymptotically optimal controls with optimality
gaps of $o(\sqrt{n})$ where $n$ is the system size, for example, the
number of servers. We construct a sequence of asymptotically optimal
controls where the optimality gap grows logarithmically with the system
size. Our analysis relies on a sequence of Brownian control problems,
whose refined structure helps us achieve the improved optimality gaps.
\end{abstract}

%
\begin{keyword}[class=AMS]
\kwd{60K25}
\kwd{90B20}
\kwd{90B36}
\kwd{49L20}
\kwd{60F17}.
\end{keyword}
\begin{keyword}
\kwd{Multiclass queues}
\kwd{many servers}
\kwd{heavy-traffic}
\kwd{Halfin--Whitt regime}
\kwd{optimal control}
\kwd{optimality gaps}
\kwd{asymptotic optimality}
\kwd{strong approximations for queues}.
\end{keyword}

\end{frontmatter}

\section{Introduction}\label{sec1}

Queueing models with many-servers are prevalent in modeling call
centers and other large-scale service systems. They are used for
optimizing staffing and making dynamic control decisions. The
complexity of the underlying queueing model renders such optimization
problems intractable for exact analysis, and one needs to resort to
approximations. A prominent mode of approximate analysis is to study
such systems in the so-called Halfin--Whitt (HW) heavy-traffic regime;
cf. \cite{HaW81}. Roughly speaking, the analysis of a queueing system
in the HW regime proceeds by scaling up the number of servers and the
arrival rate of customers in such a way that the system load approaches
one asymptotically. To be more specific, instead of considering a
single system, one considers a sequence of (closely related) queueing
systems indexed by a parameter $n$ along which the arrival rates and
the number of servers scale up so that the system traffic intensity~$\rho^n$
satisfies
%
%
\begin{equation}\label{eqHTcond}
\sqrt{n} (1-\rho^n)\rightarrow\beta\qquad\mbox{as } n\tinf.
\end{equation}

In the context of dynamic control, passing to a formal limit of the
(properly scaled) system dynamics equations as $n\tinf$ gives rise to
a \textit{limit} diffusion control problem, which is often more tractable
than the original dynamic control problem it approximates.
The approximating diffusion control problem typically provides useful
structural insights and guides the design of good policies for the
original system. Once a candidate policy is proposed for the original
problem of interest, its asymptotic performance can be studied in the
HW regime. The ultimate goal is to establish that the proposed policy
performs well. To this end, a useful criterion is the notion of
asymptotic optimality, which provides assurance that the optimality gap
associated with the proposed policy vanishes asymptotically \textit{under
diffusion scaling} as $n\tinf$. Hence, asymptotic optimality in this
context is equivalent to showing that the optimality gap is $o(\sqrt{n})$.

A central reference for our purposes is the recent paper by Atar,
Mandelbaum and Reiman \cite{AMR02}, where the authors apply all steps
of the above scheme to the important class of problems of dynamically
scheduling a multiclass queue with many identical servers in the HW
regime. Specifically, \cite{AMR02} considers a sequence of systems
indexed by the number of servers $n$, where the number of servers and
the arrival rates of the various customer classes increase with $n$
such that the heavy-traffic condition holds; cf. equation (\ref
{eqHTcond}). Following the scheme described above, the authors derive
an approximate diffusion control problem through a formal limiting
argument. They then show that the diffusion control problem admits an
optimal Markov policy, and that the corresponding HJB equation (a
semilinear elliptic PDE) has a unique classical solution. Using the
Markov control policy and the HJB equation, the authors propose
scheduling control policies for the original (sequence of) queueing
systems of interest. Finally, they prove that the proposed sequence of
policies is asymptotically optimal under diffusion scaling. Namely, the
optimality gap of the proposed policy for the $n$th system is~%
$o(\sqrt{n})$. A similar approach is applied to more general networks
in \cite{atar2005scheduling}. In this paper, we study a similar
queueing system (see Section \ref{secmodel}). Our goal, however, is
to provide an improved optimality gap which, in turn, requires a
substantially different scheme than the one alluded to above.

Approximations in the HW regime for performance analysis have been used
extensively for the study of fixed policies. Given a particular policy,
it may often be difficult to calculate various performance measures in
the original queueing system. Fortunately, the corresponding
approximations in the HW regime are often more tractable. The machinery
of strong approximations (cf. Cs\"{o}rgo and Horv\'{a}th \cite{csorgo})
often plays a central role in such analysis. In the context
of many-server heavy-traffic analysis, with strong approximations, the
arrival and service processes (under suitable assumptions on the
inter-arrival and service times) can be approximated by a diffusion
process so that the approximation error on finite intervals is $O(\log
n)$ (where $n$ is the number of servers as before). Therefore, it is
natural to expect that, under a given policy, the error in the
diffusion approximations of the various performance metrics is $O(\log
n)$, which is indeed verified for various settings in the literature
(see, e.g., \cite{MMR98}).

A natural question is then whether one can go beyond the analysis of
fixed policies and achieve an optimality gap that is logarithmic in $n$
also under dynamic control, improving upon the usual optimality gap of
$o(\sqrt{n})$. More specifically, can one propose a sequence of
policies (one for each system in the sequence) where the optimality gap
for the policy (associated with the $n$th system) is logarithmic
in~$n$? While one hopes to get logarithmic optimality gaps as suggested by
strong approximations, it is not a priori clear if this can be achieved
under dynamic control. The purpose of this paper is to provide a
resolution to this question. Namely, we study whether one can establish
such a strong notion of asymptotic optimality and if so, then how
should one go about constructing policies which are asymptotically
optimal in this stronger sense.

Our results show that such strengthened bounds on optimality gaps can
be attained. Specifically, we construct a sequence of asymptotically
optimal policies, where the optimality gap is logarithmic in $n$. Our
analysis reveals that identifying (a sequence of) candidate policies
requires a new approach. To be specific, we advance a sequence of
diffusion control problems (as opposed to just one) where the diffusion
coefficient in each system depends on the state and the control. This
is contrary to the existing work on the asymptotic analysis of queueing
systems in the HW regime. In that stream of literature, the diffusion
coefficient is typically a (deterministic) constant. Indeed, Borkar
\cite{borkar2005controlled} views the constant diffusion coefficient
as a characterizing feature of the problems stemming from the
heavy-traffic approximations in the HW regime. Interestingly, it is
essential in our work to have the diffusion coefficient depend on the
state and the control for achieving the logarithmic optimality gap. In
essence, incorporating the impact of control on the diffusion
coefficient allows us to track the policy performance in a more refined manner.

While the novelty of having the diffusion coefficient depend on the
control facilitates better system performance, it also leads to a more
complex diffusion control problem. In particular, the associated HJB
equation is fully nonlinear; it is also nonsmooth under a linear
holding cost structure. In what follows, we show that each of the HJB
equations in the sequence has a unique smooth solution on bounded
domains and that each of the diffusion control problems (when
considered up to a stopping time) admits an optimal Markov control
policy. Interpreting this solution appropriately in the context of the
original problem gives rise to a policy under which the optimality gap
is logarithmic in $n$. As in the performance analysis of fixed
policies,\vadjust{\goodbreak} strong approximations will be used in the last step, where we
propose a sequence of controls for the original queueing systems, and
show that we achieve the desired performance. However, it is important
to note that strong approximation results alone are not sufficient for
our results. Rather, for the improved optimality gaps we need the
refined properties of the solutions to the HJB equations. Specifically,
gradient estimates for the sequence of solutions to the HJB equations
(cf. Theorem \ref{thmHJB1sol}) play a central role in our proofs.

Our analysis restricts attention to a linear holding cost structure.
However, we expect the analysis to go through for some other cost
structures including convex holding costs. Indeed, the analysis of the
convex holding cost case will probably be simpler as one tends to get
``interior'' solutions in that case as opposed to the corner solutions
in the linear cost case, which causes nonsmoothness. One could also
enrich the model by allowing abandonment. We expect the analysis to go
through with no major changes in these cases as well; see the
discussion of possible extensions in Section \ref{secconclusions}.
For purposes of clarity, however, we chose not to incorporate these
additional/alternative features because we feel that the current set-up
enables us to focus on and clearly communicate the main idea: the use
of a novel Brownian model with state/control dependent diffusion
coefficient to obtain improved optimality gaps.

\subsection*{Organization of the paper} Section \ref{secmodel}
formulates the model and states the main result. Section \ref
{secresult} introduces a (sequence of) Brownian control problem(s),
which are then analyzed in Section \ref{secADCP}. A performance
analysis of our proposed policy appears in Section \ref{sectracking}.
The major building blocks of the proof are combined to establish the
main result in Section \ref{seccombining} and some concluding remarks
appear in Section~\ref{secconclusions}.

\section{Problem formulation}\label{secmodel}

\begin{figure}

\includegraphics{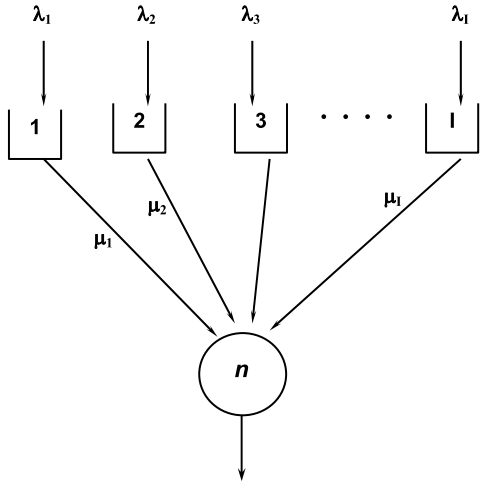}

\caption{A multiclass queue with many servers.}
\label{figv}
\end{figure}

We consider a queueing system with a single server-pool consisting of
$n$ identical servers (indexed from 1 to $n$) and a set $\I=\{1,\ldots
,I\}$ of job classes as depicted in Figure \ref{figv}. Jobs of
%
class-$i$ arrive according to a Poisson process with rate $\lambda_i$
and wait in their designated queue until their service begins. Once
admitted to service, the service time of a class-$i$ job is distributed
as an exponential random variable with rate $\mu_i>0$. All service and
interarrival times are mutually independent.

\subsection*{Heavy-traffic scaling}

We consider a sequence of systems indexed by the number of servers $n$.
The superscript $n$ will be attached to various processes and
parameters to make the dependence on $n$ explicit. (It will be omitted
from parameters and other quantities that do not change with $n$.) We
assume\vspace*{1pt} that $\lambda_i^n=a_i\lambda^n$ for all~$n$, where $\lambda
^n$ is the total arrival rate and $a_i>0$ for $i\in\I$ with $\sum
_{i}a_i=1$. This assumption is made for simplicity of notation and
presentation. Nothing changes\vspace*{1pt} in our results if one assumes, instead,
that $\lambda_i^n/n\rightarrow\lambda_i$ and $\sqrt{n}(\lambda
_i^n/n-\lambda_i)\rightarrow\hat{\lambda}_i$ as $n\tinf$ where
$\lambda_i/\sum_{k\in\I}\lambda_k=a_i>0$.\vadjust{\goodbreak}

The nominal load in the $n$th system is then given by
\[
R^n=\sum_{i}\frac{\lambda_i^n}{\mu_i}=\lambda^n\sum_{i}\frac
{a_i}{\mu_i},
\]
so that defining $\bar{\mu}=[\sum_{i}a_i/\mu_i]^{-1}$ we have that
$R^n=\lambda^n/\bar{\mu}$, which corresponds to the nominal number
of servers required to handle all the incoming jobs. The heavy-traffic
regime is then imposed by requiring that the number of servers deviates
from the nominal load by a term that is a square root of the nominal
load. Formally, we impose this by assuming that $\lambda^n$ is such
that
%
%
\begin{equation} \label{eqHTstaff}
n= R^n+\beta\sqrt{R^n}
\end{equation}
for some
$\beta\in(-\infty,\infty)$ that does not scale with $n$. Also, we
define the relative load imposed on the system by class-$i$ jobs,
denoted by $\nu_{i}$, as follows:
%
%
\begin{equation}\label{eqrhodefin}
\nu_i =\frac{a_i/\mu_i}{\sum
_{k\in\I}a_k/\mu_k}.
\end{equation}
Note that $\sum_{i\in\I}\nu_i=1$, and $\nu_i n$ can be interpreted
as a first-order (fluid) estimate for the number of servers busy
serving class-$i$ customers.

\subsection{System dynamics}

Let $Q_i\lam(t)$ and $X_i\lam(t)$ denote the number of class-$i$ jobs
in the queue and in the system, respectively, at time $t$ in the $n$th
system. Similarly, let $Z_i\lam(t)$ denote the number of servers
working on class-$i$ jobs at time $t$. Clearly, for all $i, n,
t$, the following holds:
\[
X_i\lam(t)=Z_i\lam(t)+Q_i\lam(t).\vadjust{\goodbreak}
\]

In our setting, a control corresponds to determining how many of the
\mbox{class-$i$} jobs currently in the system are placed in queue and in
service for $i\in\I$. We take the process $Z^n$ as our control in the
$n$th system. Note that one can equivalently take the queue length
process $Q^n$ as the control. (The knowledge of either process is
sufficient to pin down the evolution of the system given the arrival,
service processes and the initial conditions.) Clearly, the control
process must satisfy certain requirements for admissibility, including
the usual nonanticipativity requirement. We defer a precise
mathematical definition of admissible controls for now (see Definition
\ref{definadmissiblecontrols}). However, it should be clear that,
given the process $Z^n$, one can construct the other processes of interest.

To be specific, consider a complete probability space $(\Omega
,\mathcal{F},\mathbb{P})$ and $2I$ mutually independent
\textit{unit-rate} Poisson processes $(\mathcal{N}_i^a(\cdot), \mathcal
{N}_i^d(\cdot), i\in\I)$ on that space. Given the \textit{primitives}
$(\mathcal{N}_i^d(\cdot),\mathcal{N}_i^a(\cdot),X_i\lam(0),Z_i\lam
(0);i\in\I)$ and the control process~$Z^n$, we construct the
processes $X^n, Q^n$ as follows: for $t \geq0$ and $i \in\mathcal{I}$
%
%
\begin{eqnarray} \label{eqdynamics1}
X_i\lam(t)&=&X_i\lam(0)+\mN_i^a(\lambda_i^n t)-\mN_i^d\biggl( \mu
_i\int_0^t Z_i\lam(s)\,ds \biggr),\\
Q_i\lam(t)&=&X_i\lam(t)-Z_i\lam(t).
\end{eqnarray}
The processes $Z^n, Q^n, X^n$ must jointly satisfy the constraints
%
%
\begin{equation}\label{eqnon-negativity}
(Q\lam(t), X\lam(t),Z\lam(t))\in\mathbb{Z}_+^{3I},\qquad
e\cdot Z\lam(t)\leq n,
\end{equation}
where $e$ is the $I$-dimensional vector of ones.

Controls can be preemptive or nonpreemptive. Under a nonpreemptive
control, a job that is assigned to a server keeps the server busy until
its service is completed. In particular, given a nonpreemptive control
$Z^n$, the process $Z_i^n$ can decrease only through service
completions of class-$i$ jobs. In contrast, the class of preemptive
controls is broader. While it includes nonpreemptive policies, it also
includes controls that (occasionally) may preempt a job's service. The
preempted job is put back in the queue and its service is resumed at a
later time (possibly by a different server). Hence, the class of
preemptive controls subsumes the class of nonpreemptive ones (which is
also immediate from Definition 1 in \cite{AMR02}) and the cost of an
optimal policy among preemptive ones gives a lower bound for that among
the nonpreemptive ones.

In what follows, we will largely focus on preemptive controls, which
are easier to work with, and derive a specific policy which is near
optimal in that class. The specific policy we derive is, however,
nonpreemptive, and therefore, is near optimal among the nonpreemptive
policies as well. More specifically, the policy we propose belongs to a
class which we refer to as \textit{tracking policies.}

To facilitate the definition of tracking policies, define $\mathcal{U}
\subset\mathbb{R}_+^I$ as
%
%
\begin{equation}\label{eqmUdefin}
\mathcal{U}=\biggl\{u\in\bbR_+^I\dvtx
\sum_{i}u_i=1\biggr\}.
\end{equation}
Also, for all $i$ and $t
\geq0$, let
%
%
\begin{equation} \label{eqtildeXdefin}
\check{X}_i\lam(t)=X_i\lam(t)-\nu_in.
\end{equation}
Hence, the process $\check{X}_i^n$ captures the
oscillations of the process $X_i\lam$ around its ``fluid''
approximation $\nu_i n$. Throughout our analysis, for $x\in\bbR$ we
let $(x)^+=\max\{0,x\}$ and $(x)^-=\max\{0,-x\}$.
\begin{defin}\label{defintracking}
Given a function $h\dvtx\bbR^I\to\mathcal{U}$, an $h$-tracking policy
makes resource allocation decisions in the $n$th system as follows:
\begin{longlist}
\item It is nonpreemptive. That is, once a server starts
working on a job, it continues without interruption until that job's
service is completed.
\item It is work conserving. That is, the number of busy
servers satisfies $e\cdot Z^n(t)=(e\cdot X^n(t))\wedge n$ for all $t>
0$. In particular, no server is idle as long as there are $n$ or more
jobs in the system.
\item When a class-$i$ job arrives to the system it joins the
queue of class $i$ if all servers are busy processing other jobs.
Otherwise, the lowest-indexed idle server starts working on that job.
\item A server that finishes processing a job at a time $t$,
idles if all queues are empty. Otherwise, she starts working on a job
of class $i\in\mathcal{K}(t-)$ with probability $\lambda_i^n/\sum
_{k\in\mathcal{K}(t-)}\lambda_k^n$, where, for $t>0$, the set
$\mathcal{K}(t-)$ is defined by
%
%
\begin{equation}\label{eqmathKdefin}
\mathcal{K}(t-)=\bigl\{k\in\I\dvtx Q_k(t)-h_k(\check
{X}^n(t-))\bigl(e\cdot\check{X}^n(t-)\bigr)^+>0\bigr\}.
\end{equation}
Finally, if $(e\cdot\check{X}^n(t-))^+>0$ and $\mathcal
{K}(t-)=\varnothing$, she picks for service a customer from the lowest
index nonempty queue.
\end{longlist}
\end{defin}
\begin{rem}\label{remrandomization}
For our optimality-gap bounds and, in particular, for the proof of
Theorem \ref{thmSSC} it is important that the policy be such that
each of the job classes in the set $\mathcal{K}(t)$ gets a sufficient
share of the capacity. This prevents excessive oscillation of the
queues that may compromise the optimality gaps. Such oscillations could
arise if, for example, the policy chooses for service a~job of class
\[
i=\min\argmax_{k\in\I}\bigl\{Q_k(t-)-h_k(\check{X}^n(t-))\bigl(e\cdot
\check{X}^n(t-)\bigr)^+\dvtx Q_k^n(t-)>0\bigr\}.
\]
Randomization is just one way
to overcome such oscillations and, as the proofs (specifically that
of Theorem \ref{thmSSC}) reveal, any choice rule that guarantees a
sufficient share of the capacity to a class in $\mathcal{K}(t-)$ will
suffice.
\end{rem}

Our main result shows that a (nonpreemptive) tracking policy can
achieve a near optimal performance among preemptive policies. Note that
in our setting under\vadjust{\goodbreak} preemption, one can restrict attention to
work-conserving policies, that is, policies under which the servers
never idle as long as there are jobs to work on.\footnote{By a coupling
argument, this can be shown to hold with general queueing costs
provided that there are no abandonments and that the service times are
exponential; see, for example, the coupling argument on page 1126 of
\cite{AMR02}.} More precisely, a control is work conserving if the
following holds for all $t>0$:
%
%
\begin{equation}\label{eqworkconservation}
e\cdot Q\lam(t)=\bigl(e\cdot\check{X}\lam(t)\bigr)^+.
\end{equation}

Hereafter, we focus on work-conserving controls. Each such control can
be mapped into a ratio control, which specifies what fraction of the
total number of jobs in queue belongs to each class. To that end, let
%
%
\begin{equation}\label{eqUQmap}
U_i\lam(t)=\frac{Q_i\lam(t)}{(e\cdot Q\lam(t))\vee1} .
\end{equation}
Note that the original control $Z^n$ can be recovered from the ratio
control~$U^n$ as follows:
\[
Z_{i}^{n}(t) = X_i^n(t) - U_{i}^{n}(t) \bigl(e\cdot\check{X}^n(t)\bigr)^+ .
\]
Equations (\ref{eqdynamics1})--(\ref{eqnon-negativity}) can then be
replaced by
%
%
\begin{eqnarray}
\label{eqdynamics2}
X_i\lam(t)&=&X_i\lam(0)+\mN_i^a(\lambda_i^n t) \nonumber\\[-8pt]\\[-8pt]
&&{} -\mN_i^d\biggl(\mu_i \int_0^t \bigl(
X_i\lam(s)-U_i\lam(s) \bigl(e\cdot\check{X}^n(t)\bigr)^+ \bigr) \,ds\biggr), \nonumber\\
Q_i\lam(t)&=&U_i\lam(t)\bigl(e\cdot\check{X}^n(t)\bigr)^+, \\
Z_i\lam(t)&=&X_i\lam(t)-Q_i\lam(t),\\
\check{X}_i\lam(t)&=&X_i\lam(t)-\nu_i n,\\
\label{eqnon-negativity2}
U\lam(t)&\in&\mathcal{U},\qquad Q\lam(t)\in\mathbb{Z}_+^I,\qquad X\lam(t)\in
\mathbb{Z}_+^I.
\end{eqnarray}

Define the filtration
\[
\bar{\mathcal{F}}_t=\sigma\{\mN_i^a(s),\mN_i^d(s);i\in\I, s\leq
t\}
\]
and the $\sigma$-field
%
%
\begin{equation}\label{eqcheckFdefin}
\bar{\mathcal{F}}_{\infty}=\bigvee
_{t\geq0} \bar{\mathcal{F}}_t.
\end{equation}
Informally,
$\bar{\mathcal{F}}_{\infty}$ contains the information about the
entire evolution of the processes $(\mN_i^a,\mN_i^d,i\in\I)$. A
natural notion of admissibility requires that the control is
nonanticipative so that it only uses historical information about the
process $X^n$ and about the arrivals and service completions up to the
decision epoch. To accommodate randomized policies (as the $h$-tracking
policy) we allow the control to use other information too as long as
this information is independent of $\bar{\mathcal{F}}_{\infty}$.
\begin{defin}\label{definadmissiblecontrols}
A process $U=(U_i(t), t\geq0, i \in\mathcal{I})$
is a ratio control for the $n$th system if there exists a process
$\mathbb{X}\lam=(X\lam,Q\lam,Z\lam,\check{X}\lam)$ such that,
together with the primitives, $(\mathbb{X}\lam,U)$ satisfies
(\ref{eqdynamics2})--(\ref{eqnon-negativity2}). The process
$U$ is an admissible ratio control if, in addition, it is adapted to
the filtration $\mathcal{G}\vee\mathcal{F}_t\lam$ where
\begin{eqnarray*}
\mathcal{F}_t\lam&=&\sigma\biggl\{\mN_i^a(\lambda
_i^n s),X_i\lam(s),\mu_i\int_0^s Z_i\lam(u)\,du,\\
&&\hphantom{\sigma\biggl\{}
\mN_i^d\biggl(\mu_i\int_0^s Z_i\lam(u)\,du\biggr); i\in
\I,0\leq s\leq t\biggr\},
\end{eqnarray*}
and $\mathcal{G}$ is a $\sigma$-field that is independent of $\bar
{\mathcal{F}}_{\infty}$. The process $\mathbb{X}^n$ is then said to
be the queueing process associated with the ratio control $U$. We let
$\Pi\lam$ be the set of admissible ratio controls for the $n$th
system.
\end{defin}

Ratio controls are work conserving by definition, but they need not be
nonpreemptive in general. However, note that given a function $h\dvtx \bbR
^I\to\mathcal{U}$, the (nonpreemptive) $h$-tracking policy
corresponds to a ratio control $U_h$, which is nonpreemptive. To be
specific, given the primitives and the $h$-tracking policy, one can
construct the corresponding queueing process $\mathbb
{X}^n=(X^n,Q^n,Z^n,\check{X}^n)$ (see the construction after Lemma
\ref{lemstrongappbounds}). Then the ratio control $U_h$ is
constructed using the relation (\ref{eqUQmap}) so that $\mathbb
{X}^n$ and $U_h$ jointly satisfy (\ref{eqdynamics2})--(\ref
{eqnon-negativity2}). Hence, one can speak of the ratio control and
the queueing process associated with an $h$-tracking policy. Note that
since the tracking policy makes resource allocation decisions using
only information on the state of the system at the decision epoch
(together with a randomization that is independent of the history), the
resulting ratio control is admissible in the sense of Definition \ref
{definadmissiblecontrols}. The terms ratio control and $h$-tracking
policy appear in several places in the paper. It will be clear from the
context whether we refer to an arbitrary ratio control or to one
associated with an $h$-tracking policy.

We close this section by stating the main result of the paper. To that
end, let
%
%
\begin{equation}\label{eqmathXdefin}
\mathcal{X}^n=\{(x,q)\in\mathbb{Z}_+^{2I}\dvtx
q=u(e\cdot
x-n)^+\mbox{ for some } u\in\mathcal{U}\}.
\end{equation}
That is, $\mathcal{X}^n$ is the set on which $(X^n, Q^n)$ can
obtain values under work conservation. In this set $e\cdot q=(e\cdot
x-n)^+$ so that positive queue and idleness do not co-exist. We let
$\Ex_{x,q}^{U}[\cdot]$ denote the expectation with respect to the
initial condition $(X^n(0),Q^n(0))=(x,q)$ and an admissible ratio
control $U$. Given a ratio control $U$ and initial conditions $(x,q)$,
the expected infinite horizon discounted cost in the $n$th system
is given by
%
%
\begin{equation}\label{eqcost1}
C\lam(x,q,U)=\Ex_{x,q}^{U}\biggl[\int_0^{\infty
}e^{-\gamma s}
c\cdot Q\lam(s)\,ds\biggr],
\end{equation}
where $c=(c_1,\ldots,c_I)'$ is the strictly positive vector of holding
cost rates and $\gamma> 0$ is the discount rate. For $(x,q) \in
\mathcal{X}^n$, the value function is given by
\[
V\lam(x,q)=\inf_{U\in\Pi\lam}\Ex_{x,q}^{U}\biggl[\int_0^{\infty
} e^{-\gamma s}c\cdot Q\lam(s)\,ds\biggr].
\]

We next state our main result.
\begin{theorem}\label{thmmain} Fix a sequence $\{(x^n,q^n),n\in\bbZ_+\}
$ such that
\mbox{$(x^n,q^n)\in\mathcal{X}^n$} and $|x^n-\nu n| \leq M \sqrt{n}$ for
all $n$ and some $M>0$. Then, there exists a sequence of tracking
functions $\{h^n,n\in\bbZ_+\}$ together with constants $C,k>0$ (that
do not depend on $n$) such that
\[
C\lam(x^n,q^n,U_h^n)\leq V\lam(x^n,q^n)+C\log^{k} n \qquad\mbox{for all
} n,
\]
where $U_h^n$ is the ratio control associated with the $h^{n}$-tracking
policy.
\end{theorem}

The constant $k$ in our bound may depend on all system and cost
parameters but not on $n$. In particular, it may depend on $(\mu
_i,c_i,a_i;i\in\I)$ and $\beta$. Its value is explicitly defined
after the statement of Theorem \ref{thmHJB1sol}.

Theorem \ref{thmmain} implies, in particular, that the optimal
performance for nonpreemptive policies is close to that among the
larger family of preemptive policies. Indeed, we identify a
nonpreemptive policy (a tracking policy) in the queueing model whose
cost performance is close to the optimal value of the preemptive
control problem.

The rest of the paper is devoted to the proof of Theorem \ref{thmmain},
which proceeds by studying a sequence of auxiliary Brownian
control problems. The next subsection offers a heuristic derivation and
a justification for the relevance of the sequence of Brownian control
problems to be considered in later sections.

\subsection{Toward a Brownian control problem}
We proceed by deriving a sequence of approximating Brownian control
problems heuristically, which will be instrumental in deriving a
near-optimal policy for our original control problem. It is important
to note that we derive an approximating Brownian control problem for
each $n$ as opposed to deriving a single approximating problem (for the
entire sequence of problems). This distinction is crucial for achieving
an improved optimality gap for $n$ large because it allows us to tailor
the approximation to each element of the sequence of systems.

To this end, let
\[
l_i^n=\lambda_i^n-\mu_i\nu_i n \qquad\mbox{for } i\in\I.
\]
Fixing an admissible control $U^n$ for the $n$th system [and
centering as in (\ref{eqtildeXdefin})], we can then write (\ref
{eqdynamics2}) as
%
%
\begin{equation}\label{eqcheckXdynamics}
\check{X}_i\lam(t)=\check{X}_i\lam(0)+l_i\lam t-\mu
_i\int_0^t
\bigl( \check{X}_i\lam(s)-U^n_i(s)\bigl( e\cdot\check{X}\lam(s)\bigr)^+
\bigr) \,ds+\check{W}_i\lam(t),\hspace*{-30pt}
\end{equation}
where
%
%
\begin{eqnarray}\label{eqWtildedefin}
\check{W}_i\lam(t)&=&\mN_i^a(\lambda_i^n t)-\lambda_i^n
t+\mu_i\int_0^t \bigl( \check{X}_i\lam(s)-U^n_i(s)\bigl( e\cdot\check
{X}\lam(s)\bigr)^+ \bigr) \,ds\nonumber\\[-8pt]\\[-8pt]
&&{}-
\mN_i^d\biggl(\mu_i\int_0^t \bigl( \check{X}_i\lam(s)+\nu_i
n-U^n_i(s)\bigl( e\cdot\check{X}\lam(s)\bigr)^+ \bigr) \,ds \biggr).\nonumber
\end{eqnarray}
In words, $\check{W}_i\lam(t)$ captures the deviations of the Poisson
processes from their means. It is\vspace*{1pt} natural to expect that an
approximation result of the following form will hold: $(\check
{X}_i^n,\check{W}_i^n;i\in\I)$ can be approximated by $(\hX
_i^n,\hW_i^n;i\in\I)$ where
\begin{eqnarray*}
\hat{X}_i\lam(t)&=&\hat{X}_i\lam(0)+l_i^n t -\mu_i\int_0^t \bigl(
\hat{X}_i\lam(s)-U^n_i(s)\bigl(e\cdot\hat{X}\lam(s)\bigr)^+ \bigr) \,ds+\hat
{W}_i\lam(t),
\\
\hat{W}_i(t)&=&\tilde{B}_i^a(\lambda_i^n t)+\tilde{B}_i^S\biggl(\mu
_i\int_0^t \bigl( \hat{X}_i\lam(s)+\nu_i n-U^n_i(s)\bigl( e\cdot\hat
{X}\lam(s)\bigr)^+ \bigr) \,ds\biggr)
\end{eqnarray*}
and $\tilde{B}^a,\tilde{B}^s$ are $I$-dimensional independent
standard Brownian motions. Moreover, by a time-change argument we can
write (see, e.g., Theorem 4.6 in \cite{KaS91})
%
%
\begin{eqnarray}\label{eqhatX}
\hat{X}_i\lam(t) &=&\hat{X}_i\lam
(0)+l_i^n t -\mu_i\int_0^t \hat{X}_i\lam(s)-U^n_i(s)\bigl(e\cdot\hat
{X}\lam(s)\bigr)^+ \,ds\nonumber\\[-8pt]\\[-8pt]
&&{}+\int_0^t \sqrt{\lambda_i^n+\mu_i\bigl( \hat
{X}_i\lam(s)+\nu_i n-U^n_i(s)\bigl( e\cdot\hat{X}\lam(s)\bigr)^+
\bigr)}\,dB_i(s),\nonumber\hspace*{-30pt}
\end{eqnarray}
where $B$ is an $I$-dimensional standard Brownian motion constructed by setting
\begin{eqnarray*}B_i(t)&=& \int_0^t\frac{d\tilde{B}_i^{S}(\mu
_i\int_0^s ( \hat{X}_i\lam(u)+\nu_i n-U^n_i(u)( e\cdot\hat
{X}\lam(u))^+ ) \,du)}{\sqrt{\mu_i ( \hat
{X}_i\lam(s)+\nu_i n-U^n_i(s)( e\cdot\hat{X}\lam(s))^+
)}}\\
&&{} + \frac{\tilde{B}_i^a(\lambda_i^nt)}{\lambda_i^n t}.
\end{eqnarray*}

Taking a leap of faith and arguing heuristically, we next consider a
Brownian control problem with the system dynamics
%
%
\begin{equation}\label{eqnbmdefn1}
\hat{X}\lam(t)=x+\int_0^t b\lam(\hat{X}\lam
(s),\hat
{U}^n(s))\,ds+\int_0^t \sigma\lam(\hat{X}\lam(s),\hat
{U}^n(s))\,dB(t),\hspace*{-30pt}
\end{equation}
where $\hat{U}^n$ will be
an admissible control for the Brownian system and
%
%
\begin{equation}
\label{eqnbmdefn30}
b_i\lam(x,u)=l_i\lam-\mu_i\bigl(x_i-u_i(e\cdot x)^+\bigr)
\end{equation}
and
\begin{equation}
\label{eqnbmdefn3}
\sigma_i\lam(x,u)=\sqrt{\lambda_i^n+\mu_i \nu_in +\mu
_i\bigl(x_i-u_i(e\cdot x)^+\bigr)}.
\end{equation}
Note that the Brownian control problem will only be used to propose a~candidate policy, whose near optimality will be verified from first
principles without relying on the heuristic derivations of this section.

To repeat, the preceding definition is purely formal and provided only
as a means of motivating our approach. In what follows, we will
directly state and analyze an auxiliary Brownian control problem
motivated by the above heuristic derivation. The analysis of the
auxiliary Brownian control problem lends itself to constructing near
optimal policies for our original control problem. To be more specific,
the system dynamics equation (\ref{eqnbmdefn1}), and in particular,
the fact that its variance is state and control dependent, is crucial
to our results. Indeed, it is this feature of the auxiliary Brownian
control problems that yields an improved optimality gap.

Needless to say, one needs to take care in interpreting (\ref
{eqnbmdefn1})--(\ref{eqnbmdefn3}), which are meaningful only up
to a suitably defined hitting time. In particular, to have $\sigma^n$
well defined, we restrict attention to the process while it is within
some bounded domain. Actually, it suffices for our purposes to fix
$\kappa>0$ and $m\geq3$ and consider the Brownian control problem
only up to the hitting time of a ball of the form
%
%
\begin{equation}\label{eqBkappadefin}
\mB_{\kappa}^n=\bigl\{x\in\mathbb{R}^I\dvtx |x|< \kappa
\sqrt{n}\log^m
n\bigr\},
\end{equation}
where \mbox{$|\cdot|$} denotes the Euclidian
norm. We will fix the constant $m$ throughout and suppress the
dependence on $m$ from the notation. Setting
%
%
\begin{equation}\label{eqnkappa}
n(\kappa)=\inf\{n\in
\bbZ_+\dvtx \sigma^n(x,u)\geq1 \mbox{ for all } x\in\mB_{\kappa}^n,
u\in\mathcal{U}\},
\end{equation}
the diffusion coefficient is strictly positive for all $n \geq n(\kappa
)$ and $x\in\mB_{\kappa}^n$. Note that, for all $i\in\I$, $x\in
\mB_{\kappa}^n$ and $u\in\mathcal{U}$,
\[
(\sigma_i^n(x,u))^2\geq\lambda_i^n+\mu_i\nu_in-2\mu_i\kappa\sqrt
{n}\log^mn,
\]
so that $(\sigma_i^n(x,u))^2\geq\mu_i\nu_in/2\geq1$ for all
sufficiently large $n$ and, consequently, $n(\kappa)<\infty$.
\begin{rem} In what follows, and, in particular, through the proof of
Theorem \ref{thmmain}, the reader should note that while choosing the
size of the ball to be $\epsilon n$ (with $\epsilon$ small enough)
would suffice for the nondegeneracy of the diffusion coefficient, that
choice would be too large for our optimality gap proofs.
\end{rem}

\section{An approximating diffusion control problem (ADCP)}\label{secresult}

Motivated by the discussion in the preceding section, we define
admissible systems as follows.
\begin{defin}[(Admissible systems)]\label{definadmissiblesystemBrownian}
Fix $\kappa>0$, $n\in\bbZ
_+$ and $x\in\mathbb{R}^I$. We refer to
$\theta=(\Omega,\mathcal{F},(\mathcal{F}_t),\mathbb{P},\hat
{U},B)$ as an admissible $(\kappa,n)$-system if:
\begin{longlist}[(a)]
\item[(a)] $(\Omega,\mathcal{F},(\mathcal{F}_t),\mathbb{P})$ is
a complete filtered probability space.\vadjust{\goodbreak}
\item[(b)] $B(\cdot)$ is an $I$-dimensional standard Brownian motion
adapted to $(\mF_t)$.
\item[(c)] $\hat{U}$ is $\mathcal{U}$-valued, $\mF$-measurable
and $(\mF_t)$ progressively measurable.
\end{longlist}

The process $\hat{U}$ is said to be the control associated with
$\theta$. We also say that $\hat{X}$ is a controlled process
associated with the initial data $x$ and an admissible system $\theta$
if $\hat{X}$ is a continuous $(\mF_t)$-adapted process on $(\Omega
,\mathcal{F},\mathbb{P})$ such that, almost surely, for $t\leq\hat
{\tau}_{\kappa}^n$,
\[
\hat{X}(t)=x+\int_0^t b\lam(\hat{X}(s),\hat{U}(s))\,ds+\int_0^t
\sigma\lam(\hat{X}(s),\hat{U}(s))\,d\tilde{B}(t),
\]
where $b^n(\cdot,\cdot)$ and $\sigma^n(\cdot,\cdot)$ are as
defined in (\ref{eqnbmdefn30}) and (\ref{eqnbmdefn3}),
respectively, and
$\hat{\tau}_{\kappa}^n=\inf\{t\geq0\dvtx\hat{X}(t)\notin\mB_{\kappa
}^n\}$. Given $\kappa>0$ and $n\in\bbZ_+$, we let $\Theta(\kappa
,n)$ be the set of admissible $(\kappa,n)$-systems.
\end{defin}

The Brownian control problem then corresponds to optimally choosing an
admissible $(\kappa,n)$-system with associated control $(\hat
{U}(t),t\geq0)$ that achieves the minimal cost in the optimization problem
%
%
\begin{equation}\label{eqoptbrownian}
\hat{V}\lam(x,\kappa)=\inf_{\theta\in\Theta
(\kappa,n)}\Ex
_{x}^{\theta}\biggl[\int_0^{\hat{\tau}_{\kappa}^n} e^{-\gamma s}
\sum_{i\in\I}c_i \hat{U}_i(s)\bigl(e\cdot\hat{X}(s)\bigr)^+\,ds
\biggr],
\end{equation}
where $\Ex_x^{\theta}[\cdot]$ is the
expectation operator when the initial state is $x\in\mathbb{R}^I$ and
the admissible system $\theta$. Hereafter, we refer to (\ref
{eqoptbrownian}) as the \textit{ADCP on} $\mB_{\kappa}^n$.
The following lemma shows that the Definition \ref
{definadmissiblesystemBrownian} is not vacuous. The proof appears in
the \hyperref[app]{Appendix}.
\begin{lem} \label{lemexistenceofcontrolled}
Fix the initial state $x\in\mathbb{R}^I$, $\kappa>0$,
$n\geq n(\kappa)$ and an admissible $(\kappa,n)$-system $\theta$.
Then, there exists a unique controlled process $\hat{X}$ associated
with $x$ and $\theta$.
\end{lem}

To facilitate future analysis, note from the definition of $\hat{\tau
}_k^n$ and (\ref{eqoptbrownian}) that
%
%
\begin{equation}\label{eqvaluebound}
\hat{V}^n(x,\kappa)\leq
\frac{1}{\gamma}(e\cdot c)\kappa\sqrt{n}\log^m n.
\end{equation}
\begin{defin}[(Markov controls)]
\label{definmarkoviancontrols} We say that an admissible
$(\kappa,n)$-system $\theta=(\Omega,\mathcal{F},(\mathcal
{F}_t),\mathbb{P},\hat{U},B)$ with the associated controlled process
$\hat{X}^n$ induces a Markov control if there exists a function
$g^n(\cdot)\dvtx \mathcal{B}_{\kappa}^n \to\mathcal{U}$ such that
$\hat{U}(t)=g^n(\hat{X}^n(t))$ for $t \leq\hat{\tau}_{\kappa}^n$.
We extend the function $g^n$ to $\mathbb{R}^I$ as follows:
%
%
\begin{equation}
h^n(x)= \cases{
g^n(x), &\quad$x \in\mathcal{B}_{\kappa}^{n}$, \cr
e_1, & \quad otherwise,}
\end{equation}
where $e_1$ is the $I$-dimensional vector whose first component is $1$
while the others are $0$. We refer to $h^n(\cdot)$ as the tracking
function associated with the admissible system $\theta$.
\end{defin}

In what follows, a policy $\hat{U}$ will be called optimal for the
approximating diffusion control problem (ADCP) on $\mB_{\kappa}^n$ if
there exists an admissible $(\kappa,n)$-system $\theta=(\Omega
,\mathcal{F},(\mathcal{F}_t),\mathbb{P},\hat{U},B)$ such that
\[
\hat{V}\lam(x,\kappa)=\Ex_{x}^{\theta}\biggl[\int_0^{\hat{\tau
}_{\kappa}^n} e^{-\gamma s} \sum_{i\in\I}c_i \hat{U}_i(s)\bigl(e\cdot
\hat{X}(s)\bigr)^+\,ds\biggr].
\]

Recall that $X$ and $U$ are used to denote performance relevant
stochastic processes in both the Brownian model and the original
queueing model, and that we add a hat, that is, we use $\hat{X}$ and
$\hat{U}$ in the context of the Brownian model. To avoid confusion,
the reader should keep in mind that hat-processes correspond to the
ADCP while the ones with no hats correspond to the original queueing model.

\subsection*{Roadmap for the remainder of the paper} The main result in
Theorem~\ref{thmmain} builds on the following steps:
\begin{enumerate}
\item In Section \ref{secADCP} we show that for each $n$, the HJB
equation associated with the ADCP has a unique\vadjust{\goodbreak} and sufficiently smooth
solution. Using that solution we advance an optimal Markov control for
the ADCP together with the corresponding tracking function. We also
identify useful gradient bounds on the solutions to the sequence of HJB
equations; cf. Theorem~\ref{thmHJB1sol}.

\item In Section \ref{sectracking} we conduct a performance analysis
of $h$-tracking policies in the queueing system; cf. Theorem \ref{thmSSC}.

\item The result of Theorem \ref{thmSSC} together with the gradient
estimates in Theorem \ref{thmHJB1sol} are combined in a
Taylor expansion-type argument in Section \ref{seccombining} to
complete the proof of Theorem \ref{thmmain}.
\end{enumerate}

As a convention, throughout the paper we use the capital letter $C$ to
denote a constant that does not depend on $n$. The value of $C$ may
change from line to line within the proofs but it will be clear from
the context.

\section{Solution to the ADCP} \label{secADCP} This section provides
a solution for the ADCP on $\mB_{\kappa}^n$ for each $n\in\bbZ$ and
$\kappa>0$. The HJB equation is a fully nonlinear and nonsmooth PDE.
As such, it requires extra care when compared with the usual semilinear
PDEs that arise in the analysis of \textit{asymptotically} optimal
controls in the Halfin--Whitt regime. We will build on existing results
in the theory of PDEs and proceed through the following steps: (a)
establish the existence and uniqueness of classical solutions; (b)
relate this unique solution to the value function of the ADCP and (c)
establish useful gradient estimates on the solution for the HJB
equation. The last step is not necessary for existence and uniqueness
but is important for the analysis of optimality gaps.

In what follows, we fix $\kappa>0$ and $n\geq n(\kappa)$ and suppress
the dependence of the solution to the HJB equation on $n$ and $\kappa
$. The following notation is needed to introduce the HJB equation.
Given a twice continuously differentiable function $\phi$, define
\[
\phi_i=\frac{\partial\phi}{\partial x_i}\quad\mbox{and}\quad \phi_{ii}=
\frac{\partial^2 \phi}{\partial x_i^2} .
\]
Also, define the operator $A^n_u$ for $u\in\mathcal{U}$ as follows:
%
%
\begin{equation}\label{eqgendefin}
A_{u}\lam\phi= \sum_{i\in\I} b_i\lam(\cdot
,u)\phi_i+\frac
{1}{2}\sum_{i\in\I} (\sigma_i\lam(\cdot,u))^2 \phi_{ii}.
\end{equation}
Defining
\[
L(x,u)=\sum_{i\in\I}c_iu_i (e\cdot x)^+
\]
for $x\in\bbR_+^I$ and $u\in\mathcal{U}$, the HJB equation is given by
%
%
\begin{equation} \label{eqHJB0}
0=\inf_{u \in\mathcal{U}}\{
L(x,u)+A_{u}\lam\phi
(x)-\gamma\phi(x)\}.
\end{equation}
Substituting $b\lam(\cdot,\cdot)$ and $\sigma\lam(\cdot
,\cdot)$ into (\ref{eqHJB0}) gives
%
%
\begin{eqnarray} \label{eqHJB1simp}
0&=& -\gamma\phi(x) + (e\cdot x)^+\cdot\min_{i\in\I}\biggl\{
c_i+\mu_i\phi_i(x)-\frac{1}{2}\mu_i\phi_{ii}(x)\biggr\}
\nonumber\\[-8pt]\\[-8pt]
&&{}+\sum_{i\in\I} (l_i\lam-\mu_ix_i)\phi_i(x)
+\frac{1}{2}\sum_{i\in\I} \bigl(\lambda_i^n+\mu_i(\nu
_in+x_i)\bigr)\phi_{ii}(x).\nonumber
\end{eqnarray}

Our analysis of the HJB equation (\ref{eqHJB1simp}) draws on existing
results on fully nonlinear PDEs, and, in particular, the results on
Bellman--Pucci type equations; cf. Chapter 17 of \cite{TandG}.

In what follows, fixing a set $\mB\subseteq\mathbb{R}_+^I$, $\mC
^2(\mB)$ denotes the space of twice continuously differentiable
functions from $\mB$ to $\mathbb{R}$. For $u\in\mC^{2}(\mB)$, we
let~$Du$ and $D^2u$ denote the gradient and the Hessian of $u$,
respectively. The space $\mC^{2,\alpha}(\mB)$ is then the subspace
of $\mC^{2}(\mB)$ members of which also have second derivatives that
are H\"{o}lder continuous of order $\alpha$. That is, a twice
continuously differentiable function $u\dvtx \mathbb{R}^I\to\mathbb{R}$
is in $\mC^{2,\alpha}(\mB)$ if
\[
\sup_{x,y\in\mB, x\neq y} \frac{|D^2u(x)-D^2u(y)|}{|x-y|^{\alpha
}}<\infty,
\]
where \mbox{$|\cdot|$} denotes the Euclidian norm. We define
$d_x\,{=}\,\operatorname{dist}(x,\partial\mB)\,{=}\,\inf\{|x\,{-}\,y|,\allowbreak
y\,{\in}\,\partial\mB\}$ where
$\partial\mB$ stands for the boundary of $\mB$ and we let
$d_{x,z}\,{=}\,\min\{d_x,d_z\}$. Also, we define
%
%
\begin{equation}\label{equstardefin}
|u|^*_{2,\alpha,\mB}=\sum_{j=0}^2 [u]_{j,\mB}^*+\sup_{x,y\in\mB,x\neq
y}d_{x,y}^{2+\alpha}\frac{|D^2u(x)-D^2 u(y)|}{|x-y|^{\alpha}},
\end{equation}
where $[u]_{j,\mB}^*=\sup_{x\in\mB}d_x^j
|D^j u(x)|$ for $j=0,1,2$. Note that $d_x^j$ denote the $j$th power
of $d_x$ and, similarly, $d_{x,y}^{2+\alpha}$ is the $(2+\alpha)$th
power of $d_{x,y}$. Finally, we let $|u|^*_{0,\mB}=[u]_{0,\mB
}^*=\sup_{x\in\mB}|u(x)|$.

In the statement of the following theorem, $e_j$ is the $I$-dimensional
vector with $1$ in the $j$th place and zeros elsewhere. Also, $\mB
_{\kappa}^n$, $m$ and $n(\kappa)$ are as defined in (\ref
{eqBkappadefin}) and (\ref{eqnkappa}), respectively.
\begin{theorem}\label{thmHJB1sol}
Fix $\kappa>0$ and $n\geq n(\kappa)$. Then, there
exists $0<\alpha\leq1$ (that does not depend on $n$) and a unique
classical solution $\phi_{\kappa}\lam\in\mC^{0,1}(\bar{\mB}
_{\kappa}^n)\cap\mC^{2,\alpha}(\mB_{\kappa}^n)$ to the HJB
equation (\ref{eqHJB1simp}) on $\mB_{\kappa}^n$ with the
boundary condition $\phi_{\kappa}\lam=0$ on $\partial\mB_{\kappa
}^n$. Furthermore, there exists a constant $C>0$ (that does not depend
on $n$) such that
%
%
\begin{equation}\label{eqgradients0}
|\phi_{\kappa}\lam|^*_{2,\alpha,\mB_{\kappa
}^n}\leq C\sqrt
{n}\log^{k_0} n,
\end{equation}
where $k_0=4m(1+1/\alpha)$. In turn, for any $\vartheta<1$,
%
%
\begin{equation}\label{eqgradients1}
\sup_{x\in\mB_{\vartheta\kappa}^n}|D\phi_{\kappa
}^n(x)|\leq
\frac{C}{1-\vartheta}\log^{k_1} n \quad\mbox{and}\quad\sup_{x\in\mB
_{\vartheta\kappa}^n}|D^2\phi_{\kappa}^n(x)|\leq\frac
{C}{1-\vartheta}\frac{\log^{k_2}n} {\sqrt{n}}\hspace*{-26pt}
\end{equation}
with $k_1=k_0-m$ and $k_2=k_0-2m$. Also,
%
%
\begin{eqnarray}\label{eqgenbound}
&&\sup_{u\in\mathcal{U}}\biggl|\sum_{i\in\I}\bigl((\phi_{\kappa
}^n)_{ii}(y)-(\phi_{\kappa}^n)_{ii}(x)\bigr)(\sigma_i^n(x,u))^2
\biggr|\nonumber\\[-8pt]\\[-8pt]
&&\qquad\leq\frac{C}{1-\vartheta}\log^{k_1}n\nonumber
\end{eqnarray}
for all $x,y\in\mB_{\vartheta\kappa}^n$ with $|x-y|\leq1$.
\end{theorem}

Note that (\ref{eqgradients1}) follows immediately from (\ref
{eqgradients0}) through the definition of the operation \mbox{$|\cdot
|^*_{2,\alpha,\mB_{\kappa}^n}$} in (\ref{equstardefin}).
Henceforth, we will use $k_i,i=0,1,2$ for the values given in the
statement of Theorem \ref{thmHJB1sol}. Moreover, the constant $k$
appearing in the statement of Theorem \ref{thmmain} is equal to $k_0+3$.

Theorem \ref{thmHJB1sol} facilitates a verification result, which we
state next followed by the proof of Theorem \ref{thmHJB1sol}. Below,
$\hat{V}^n(x,\kappa)$ is the value function of the ADCP; cf. equation
(\ref{eqoptbrownian}).
\begin{theorem}\label{thmBrownianverification}
Fix $\kappa>0$ and $n\geq n(\kappa)$. Let $\phi
_{\kappa}^n$ be the unique solution to the HJB equation (\ref
{eqHJB1simp}) on $\mB_{\kappa}^n$ with the boundary condition $\phi
_{\kappa}\lam=0$ on $\partial\mB_{\kappa}^n$. Then,
$\phi_{\kappa}^n(x)=\hat{V}^n(x,\kappa)$ for all $x\in\mB_{\kappa
}^n$. Moreover, there exists a Markov control which is optimal for the
ADCP on $\mB_{\kappa}^n$. The tracking function $h_{\kappa}^{*,n}$
associated with this optimal Markov control is defined by $h_{\kappa
}^{*,n}(x)=e_{i^n(x)}$, where
%
%
\begin{equation}\label{eqixdefin}
i^n(x)=\min\mathop{\argmin}_{i\in\I} \biggl\{ \biggl(c_i+\mu
_i(\phi
_{\kappa}\lam)_i(x)-\frac{1}{2}\mu_i(\phi_{\kappa}\lam
)_{ii}(x)\biggr)(e\cdot x)^+\biggr\}.\hspace*{-25pt}
\end{equation}
\end{theorem}

The HJB equation (\ref{eqHJB1simp}) has two sources of
nondifferentiability. The first source is the minimum operation and the
second is the nondifferentiability of the term $(e\cdot x)^+$. The
first source of nondifferentiability is covered almost entirely by the
results in \cite{TandG}. To deal with the nondifferentiability of the
function $(e\cdot x)^+$, we use a construction by approximations. The
proof of existence and uniqueness in Theorem \ref{thmHJB1sol} follows
an approximation scheme where one replaces the nonsmooth function
$(e\cdot x)^+$ by a smooth (parameterized by $a$) function $f_a(e\cdot
x)$. We show that the resulting ``perturbed'' PDE has a unique
classical solution and that as $a\tinf$ the corresponding sequence of
solutions converges, in an appropriate sense, to a solution to~(\ref
{eqHJB1simp}) which will be shown to be unique. Note that this argument
is repeated for each fixed $n$ and $\kappa$.

To that end, given $a>0$, define
%
%
\begin{equation}\label{eqfdefin}
f_a(y)=\cases{
y, &\quad$\displaystyle y\geq\frac{1}{4a}$,\vspace*{2pt}\cr
\displaystyle ay^2+\frac{1}{2}y+\frac{1}{16a}, &\quad$\displaystyle -\frac{1}{4a}\leq y\leq\frac
{1}{4a}$,\vspace*{2pt}\cr
0, & \quad otherwise.}
\end{equation}
Replacing $(e\cdot x)^+$ with $f_a(e\cdot x)$ in (\ref{eqHJB1simp})
gives the following equation:
%
%
\begin{eqnarray} \label{eqHJB2}
0&=& - \gamma\phi(x)+f_a(e\cdot x)\cdot\min_{i\in
\I} \biggl\{ c_i+\mu_i\phi_i(x)-\frac{1}{2}\mu_i\phi_{ii}(x)
\biggr\} \nonumber\\[-8pt]\\[-8pt]
&&{}+\sum_{i\in\I} (l_i\lam-\mu_ix_i)\phi_i(x)
+\frac{1}{2}\sum_{i\in\I} \bigl(\lambda_i^n+\mu_i(\nu
_in+x_i)\bigr)\phi_{ii}(x).\nonumber
\end{eqnarray}

To simplify this further, let $\Gamma= \mB_{\kappa}^n\times\mathbb
{R}_+\times\mathbb{R}^I\times\mathbb{R}^{I\times I}$ and for all
$y\in\Gamma$, define the function
%
%
\begin{equation}\label{eqFForm}
F_a^k[y]=\min\{F_a^1[y],\ldots, F_a^I[y]\},
\end{equation}
where for $k\in\I$ and $y=(x,z,p,r)\in\Gamma$,
%
%
\begin{eqnarray}\label{eqFdefin}
F^k_a[y]&=&f_a(e\cdot x)
\biggl[c_k+\mu_kp_k-\frac{1}{2}\mu_kr_{kk}\biggr]+\sum_{i\in\I
}(l_i\lam-\mu_ix_i)p_i
\nonumber\\[-8pt]\\[-8pt]
&&{}+\frac{1}{2}\sum_{i\in\I} \bigl(\lambda_i^n+\mu_i(\nu
_in+x_i)\bigr)r_{ii}-\gamma z.\nonumber
\end{eqnarray}
Then, (\ref{eqHJB2}) can be rewritten as
%
%
\begin{equation}\label{eqFFormPDe}
F_a[x,u(x),Du(x),D^2u(x)]=0.
\end{equation}

In the following statement we use the gradient notation introduced at
the beginning of this section.
\begin{prop}\label{propsolPHJB} Fix $\kappa>0$, $n\geq n(\kappa)$ and $a>0$. A unique
classical solution $\phi_{\kappa,a}^n\in\mC^{0,1}(\bar{\mB}
_{\kappa}^n)\cap\mC^{2,\alpha}(\mB_{\kappa}^n)$ exists for
the PDE (\ref{eqHJB2}) on $\mB_{\kappa}^n$ with the boundary
condition $\phi_{\kappa,a}^n=0$ on $\partial\mB_{\kappa}^n$. Moreover,
%
%
\begin{equation}\label{eqgradients}
|\phi^n_{\kappa,a}|^*_{2,\alpha,\mB_{\kappa
}^n}\leq C |\phi
^n_{\kappa,a}|^*_{0,\mB_{\kappa}^n}\log^{k_0}n \leq\tilde{C}
\end{equation}
for $k_0=4m(1+1/\alpha)$ where
$0<\alpha\leq1$ and $C>0$ do not depend on $a$ and $n$ and $\tilde
{C}$ does not depend on $a$. Also, $\phi^n_{\kappa,a}$ is Lipschitz
continuous on the closure $\bar{\mB}_{\kappa}^n$ with a
Lipschitz constant that does not depend on $a$ (but can depend on
$\kappa$ and $n$).
\end{prop}

We postpone the proof of Proposition \ref{propsolPHJB} to the \hyperref[app]{Appendix}
and use it to complete the proof of Theorem \ref{thmHJB1sol},
followed by the proof of Theorem \ref{thmBrownianverification}.

\subsection*{\texorpdfstring{Proof of Theorem \protect\ref{thmHJB1sol}}{Proof of Theorem 4.1}}

Since we fix $n$ and $\kappa$, they will be suppressed below. We
proceed to show the existence by an approximation argument. To that
end, fix a sequence $\{a^k;k\in\bbZ\}$ with $a^k\tinf$ as $k\tinf$
and let $\phi_{a^k}$ be the unique solution to (\ref{eqHJB2}) as
given by Proposition \ref{propsolPHJB}. The next step is to show
that $\phi_{a^k}$ has a subsequence that converges in an appropriate
sense to a function $\phi$, which is, in fact, a solution to the HJB
equation (\ref{eqHJB1simp}). To that end, let
%
%
\begin{equation}\label{eqCstar}
\mC_{*}^{2,\alpha
}(\mB)=\{u\in\mC^{2,\alpha}(\mB)\dvtx |u|_{2,\alpha,\mB}^*<\infty\}.
\end{equation}
Then, $\mC_{*}^{2,\alpha}(\mB)$ is a Banach
space (see, e.g., Exercise 5.2 in \cite{TandG}). Since the bound in
(\ref{eqgradients}) is independent of $a$, we have that $\{\phi
_{a^k}\}$ is a bounded sequence in $C_{*}^{2,\alpha}(\mB)$ and hence,
contains a convergent subsequence. Let $u$ be a limit point of the
sequence $\{\phi_{a^k}\}$. Since the gradient estimates in Proposition
\ref{propsolPHJB} are independent of $a$, they hold also for the
limit function $u$, that is,
%
%
\begin{equation}\label{equ-bound}
|u|_{2,\alpha,\mB}^*\leq C|u|_{0,\mB
}^*\log^{k_0} n\leq\tilde{C}
\end{equation}
for constants
$\alpha$ and $C$ that are independent of $n$. Proposition \ref
{propsolPHJB} also guarantees that the global Lipschitz constant is
independent of $a$ so that we may conclude that $u\in\mC
^{0,1}(\bar{\mB})$ and that $u=0$ on $\partial\mB$.

We will now show that $u$ solves (\ref{eqHJB1simp}) uniquely. To show
that $u$ solves (\ref{eqHJB1simp}), we need to show that $F[u]=0$
(where $F[\cdot]$ is defined similar to $F_a[\cdot]$ with $(e\cdot
x)^+$ replacing $f_a(e\cdot x)$). To that end, let $\{a^k,k\in\bbZ\}
$ be the corresponding convergent subsequence [i.e., such that $\phi
_{a^k}\rightarrow u$ in $\mC_*^{2,\alpha}(\mB)$]. Henceforth, to
simplify notation, we write
\[
F_{\akl}[\phi_{\akl}(x)]=F_{\akl}[x,\phi_{\akl}(x),D\phi_{\akl
}(x),D^2\phi_{\akl}(x)]
\]
(and similarly for $F[\cdot]$). Fix $\delta\,{>}\,0$ and let $\mB(\delta
)\,{=}\,\{x\,{\in}\,\bbR^I\dvtx |x|\,{<}\,\kappa\sqrt{n}\log^m n\,{-}\,\delta\}$. Note that
since $\phi_{\akl}\rightarrow u$ in $\mC_*^{2,\alpha}(\mB)$ we
have, in particular, the convergence of $(\phi_{\akl}(x),D\phi_{\akl
}(x),D^2\phi_{\akl}(x))\rightarrow(u(x),Du(x),D^2u(x))$ uniformly in
$x\in\mB(\delta)$. The equicontinuity of the function
$F^a[\cdot]$\vadjust{\goodbreak}
on $\Gamma$ guarantees then that
%
%
\begin{equation}\label{eqinterim11}
|F_{\akl}[\phi_{\akl}(x)]-F_{\akl}[u(x)]|\leq
\epsilon
\end{equation}
for all $l$ large enough and $x\in\mB(\delta)$.
Note that $\sup_{x\in\bbR^I}|f_{a}(e\cdot x)-(e\cdot x)^+|\leq
\epsilon$ for all $a$ large enough so that,
%
%
\begin{equation}\label{eqinterim12}
\sup_{x\in\mB}|F_{\akl}[u(x)]-F[u(x)]|\leq\epsilon
\end{equation}
for all $l$ large enough. Combining (\ref{eqinterim11}) and (\ref
{eqinterim12}), we then have
\[
\sup_{x\in\mB}|F_{\akl}[\phi_{\akl}(x)]-F[u(x)]|\leq2\epsilon
\]
for all $l$ large enough and $x\in\mB(\delta)$. By definition
$F^{a^k}[\phi_{\akl}(x)]=0$ for all $x\in\mB$ and since $\epsilon$
was arbitrary we have that $F[u(x)]=0$ for all $x\in\mB(\delta)$.
Finally, since $\delta$ was arbitrary we have that $F[u(x)]=0$ for all
$x\in\mB$. We already argued that $u=0$ on $\partial\mB$, so that
$u$ solves (\ref{eqHJB1simp}) on~$\mB$ with $u=0$ on~$\partial\mB
$. This concludes the proof of existence of a solution to (\ref
{eqHJB1simp}) that satisfies the gradient estimates (\ref{eqgradients0}).

Finally, the uniqueness of the solution to (\ref{eqHJB1simp}) follows
from Corollary 17.2 in \cite{TandG} noting that the function
$F[x,z,p,r]$ is indeed continuously differentiable in the $(z,p,r)$
arguments and it is decreasing in $z$ for all $(x,p,r)$.

Using Theorem \ref{thmBrownianverification} [which only uses the
existence and uniqueness of the solution $\phi_{\kappa}^n(x)$ that we
already established] together with (\ref{eqvaluebound}) we have that
\[
|\phi_{\kappa}^n|_{0,\mB_{\kappa}^n}=\sup_{x\in\mB_{\kappa
}^n}\hat{V}^n(x,\kappa)\leq\frac{1}{\gamma}\kappa\sqrt{n}\log^m n.
\]
The bounds (\ref{eqgradients0}) and (\ref{eqgradients1}) now follow
from (\ref{equ-bound}) and we turn to prove (\ref{eqgenbound}).

To that end, since $\phi_{\kappa}^n$ solves
(\ref{eqHJB1simp}), fixing $x,y\in\mB_{\kappa}^n$ we have
%
%
\begin{eqnarray}\label{eqinternational}\quad
&&\biggl|\frac{1}{2}\sum_{i\in\I}
\bigl(\lambda_i^n+\mu_i(\nu_in+x_i)\bigr)(\phi_{\kappa}^n)_{ii}(x)-
\frac{1}{2}\sum_{i\in\I} \bigl(\lambda_i^n+\mu_i(\nu
_in+y_i)\bigr)(\phi_{\kappa}^n)_{ii}(y)\biggr|\nonumber\\
&&\qquad \leq\gamma|\phi_{\kappa}^n(x)-\phi_{\kappa
}^n(y)|\nonumber\\[-8pt]\\[-8pt]
&&\qquad\quad{} +\biggl|(e\cdot x)^+\cdot\min
_{i\in\I}\biggl\{ c_i+\mu_i(\phi_{\kappa}^n)_i(x)-\frac{1}{2}\mu
_i(\phi_{\kappa}^n)_{ii}(x)\biggr\}\nonumber\\
&&\qquad\quad\hphantom{{}+\biggl|}{} -(e\cdot y)^+\cdot\min_{i\in\I}\biggl\{
c_i+\mu_i(\phi_{\kappa}^n)_i(y)-\frac{1}{2}\mu_i(\phi_{\kappa
}^n)_{ii}(y)\biggr\}\biggr|.\nonumber
\end{eqnarray}
We will now bound each of the elements on the right-hand side. To that
end, let $i(x)$ be as defined in (\ref{eqixdefin}) and for each
$x,z\in\mB_{\vartheta\kappa}^n$ define
\[
M_{i(x)}^n(z)=c_{i(x)}+\mu_{i(x)}(\phi_{\kappa}^n)_{i(x)}(z)-\tfrac
{1}{2}\mu_{i(x)}(\phi_{\kappa}^n)_{i(x)i(x)}(z).
\]
Using (\ref{eqgradients1}), we have by the mean value theorem that
%
%
\begin{equation}\label{eqinterinter}
|\phi_{\kappa}^n(x)-\phi_{\kappa}^n(y)|\leq
{|x-y|\max_{i\in\I
}\sup_{z\in\mB_{\vartheta\kappa}^n}}|(\phi_{\kappa}^n)_i(z)|\leq
C\log^{k_1} n\vadjust{\goodbreak}
\end{equation}
for all $x,y\in\mB
_{\vartheta\kappa}^n$ with $|x-y|\leq1$, and we turn to bound the
second element on the right-hand side of (\ref{eqinternational}).
Here, there are two cases to consider. Suppose first that
$i(x)=i(y)=i$. Then, using (\ref{eqgradients1}) and the mean value
theorem we have
\[
|(\phi_{\kappa}^n)_{i}(x)-(\phi_{\kappa}^n)_{i}(y)|\leq{|x-y|\max
_{i\in\I}\sup_{z\in\mB_{\vartheta\kappa}^n}}|(\phi_{\kappa
}^n)_{ii}(z)|\leq C\frac{
\log^{k_2} n}{\sqrt{n}}
\]
and, in turn, that
%
%
\begin{equation} \label{eqMbound}
|M_{i}^n(x)-M_{i}^n(y)|\leq C\frac{\log^{k_2}
n}{\sqrt{n}}
\end{equation}
for all $x,y\in\mB_{\vartheta\kappa}^n$ with $|x-y|\leq1$. Now,
$|x|\vee|y|\leq\kappa\sqrt{n}\log^m n$ for all $x,y\in\mB
_{\vartheta\kappa}^n$ and, by (\ref{eqgradients1}), $\sup_{z\in
\mB_{\vartheta\kappa}^n}|(\phi_{\kappa}^n)_{ii}(z)|\vee|(\phi
_{\kappa}^n)_{i}(z)|\leq C\log^{k_1}n$ so that
%
%
\begin{eqnarray}\label{eqinterim222222}
&& |(e\cdot x)^+ M_i^n(x)-(e\cdot y)^+
M_i^n(y)|\nonumber\\
&&\qquad
\leq
\kappa\sqrt{n}\log^m n | M_i^n(x)-M_i^n(y)|+\sup_{z\in
\mB_{\vartheta\kappa}^n}
|M_i^n(z)|\\
&&\qquad\leq C\log^{k_1}n.\nonumber
\end{eqnarray}
If, on the other hand, $i(x)\neq i(y)$ then by the definition
of $i(\cdot)$,
\begin{eqnarray*}
&&c_{i(x)}+\mu_{i(x)}(\phi_{\kappa
}^n)_{i(x)}(x)-\tfrac{1}{2}\mu_{i(x)}(\phi_{\kappa
}^n)_{i(x)i(x)}(x)\\
&&\qquad \leq
c_{i(y)}+\mu_{i(y)}(\phi_{\kappa}^n)_{i(y)}(x)-\tfrac{1}{2}\mu
_{i(y)}(\phi_{\kappa}^n)_{i(y)i(y)}(x)
\end{eqnarray*}
and
\begin{eqnarray*}
&& c_{i(y)}+\mu_{i(y)}(\phi_{\kappa
}^n)_{i(y)}(y)-\tfrac{1}{2}\mu_{i(y)}(\phi_{\kappa
}^n)_{i(y)i(y)}(y)\\
&&\qquad \leq
c_{i(x)}+\mu_{i(x)}(\phi_{\kappa}^n)_{i(x)}(y)-\tfrac{1}{2}\mu
_{i(x)}(\phi_{\kappa}^n)_{i(x)i(x)}(y).
\end{eqnarray*}
That is,
%
%
\begin{equation} \label{eqyetonemore1}
M_{i(x)}^n(x)\leq M_{i(y)}^n(x) \quad\mbox{and}\quad
M_{i(y)}^n(y)\leq M_{i(x)}^n(y).
\end{equation}
Using
(\ref{eqgradients1}) as before we have for $x,y\in\mB_{\vartheta
\kappa}^n$ with $|x-y|\leq1$ and $i(x)\neq i(y)$ that
\[
\bigl| M_{i(x)}^n(x)-M_{i(x)}^n(y)\bigr|+\bigl|
M_{i(y)}^n(x)-M_{i(y)}^n(y)\bigr| \leq C\frac{\log^{k_2} n}{\sqrt{n}}.
\]
By (\ref{eqyetonemore1}) we then have that
\begin{eqnarray*}
\bigl| M_{i(x)}^n(x)-M_{i(y)}^n(y)\bigr|&\leq&\bigl|
M_{i(x)}^n(x)-M_{i(x)}^n(y)\bigr|\\
&&{}+ \bigl|
M_{i(y)}^n(x)-M_{i(y)}^n(y)\bigr| \\
&\leq& C\frac{\log^{k_2}n}{\sqrt{n}}
\end{eqnarray*}
for all such $x$ and $y$. In turn, since $|x|\vee|y|\leq\kappa\sqrt
{n}\log^m n$,
%
%
\begin{equation}\label{eqinterim222224}
\bigl| (e\cdot x)^+M_{i(x)}^n(x)- (e\cdot
y)^+M_{i(y)}^n(y)\bigr|\leq C\log^{k_1} n
\end{equation}
for $x,y\in\mB_{\vartheta\kappa}^n$ with $|x-y|\leq1$ and
$i(x)\neq i(y)$. Plugging (\ref{eqinterinter}), (\ref{eqinterim222222})
and (\ref{eqinterim222224}) into the right-hand
side of (\ref{eqinternational}) we get
%
%
\begin{eqnarray}\label{eqtheinternational}\qquad
&&\biggl|\frac{1}{2}\sum_{i\in\I} \bigl(\lambda
_i^n+\mu_i(\nu_in+x_i)\bigr)(\phi_{\kappa}^n)_{ii}(x)-
\frac{1}{2}\sum_{i\in\I} \bigl(\lambda_i^n+\mu_i(\nu
_in+y_i)\bigr)(\phi_{\kappa}^n)_{ii}(y)\biggr|\nonumber\\[-8pt]\\[-8pt]
&&\qquad
\leq C \log^{k_1} n\nonumber
\end{eqnarray}
for all $x,y\in\mB_{\vartheta\kappa}^n$ with $|x-y|\leq1$.
Finally, recall that
\[
\sigma_i\lam(x,u)=\sqrt{\lambda_i^n+\mu_i \nu_in +\mu
_i\bigl(x_i-u_i(e\cdot x)^+\bigr)}
\]
so that for all $u\in\mathcal{U}$,
\begin{eqnarray*}
\hspace*{-4pt}&&\biggl|\sum_{i\in\I}\bigl((\phi_{\kappa
}^n)_{ii}(y)-(\phi_{\kappa}^n)_{ii}(x)\bigr)(\sigma_i^n(x,u))^2\biggr|
\\
\hspace*{-4pt}&&\qquad =\biggl|\sum_{i\in\I}(\phi_{\kappa}^n)_{ii}(y)\bigl(\lambda
_i^n+\mu_i \nu_in +\mu_i\bigl(x_i-u_i(e\cdot x)^+\bigr)\bigr)\\
\hspace*{-4pt}&&\qquad\quad\hspace*{3.5pt}{}
-(\phi_{\kappa}^n)_{ii}(x)\bigl(\lambda_i^n+\mu_i \nu_in +\mu
_i\bigl(x_i-u_i(e\cdot x)^+\bigr)\bigr)\biggr|
\\
\hspace*{-4pt}&&\qquad\leq\biggl|\frac{1}{2}\sum_{i\in\I} \bigl(\lambda_i^n+\mu
_i(\nu_in+x_i)\bigr)(\phi_{\kappa}^n)_{ii}(x)\\
\hspace*{-4pt}&&\qquad\quad\hspace*{2pt}{}-
\frac{1}{2}\sum_{i\in\I} \bigl(\lambda_i^n+\mu_i(\nu
_in+y_i)\bigr)(\phi_{\kappa}^n)_{ii}(y)\biggr|\\
\hspace*{-4pt}&&\qquad\quad{} +\biggl|\frac{1}{2}\sum_{i\in\I} (\phi_{\kappa
}^n)_{ii}(x)\mu_iu_i(e\cdot x)^+-(\phi_{\kappa}^n)_{ii}(y)\mu
_iu_i(e\cdot y)^+\biggr|\\
\hspace*{-4pt}&&\qquad\quad{} +\biggl|\frac{1}{2}\sum_{i\in
\I} (\phi_{\kappa}^n)_{ii}(y)\mu_i(x_i-y_i)\biggr|.
\end{eqnarray*}
The last two terms above are bounded by $C\log^{k_1}n$ by (\ref
{eqgradients1}) and using $|x|\vee|y|\leq\kappa\sqrt{n}\log^mn$.
Together with (\ref{eqtheinternational}) this establishes (\ref
{eqgenbound}) and concludes the proof of the theorem.

\subsection*{\texorpdfstring{Proof of Theorem \protect\ref{thmBrownianverification}}{Proof of Theorem 4.2}}
Fix an initial condition $x\in\mB_{\kappa}^n$ and an admissible
$(\kappa,n)$-system $\theta=(\Omega,\mathcal{F},(\mathcal
{F}_t),\mathbb{P},\hat{U},B)$ and let $\hat{X}^n$ be the associated
controlled process.\vadjust{\goodbreak} Using It\^{o}'s lemma for the function $\varphi
(t,x)=e^{-\gamma t} \phi_{\kappa}\lam(x)$ in conjunction with the inequality
\[
L(x,u)+A_{u}\phi_{\kappa}\lam(x)-\gamma\phi_{\kappa}\lam(x)\geq
0 \qquad\mbox{for all } x\in\mB_{\kappa}^n, u\in\mathcal{U}
\]
[recall that $\phi_{\kappa}^n$ solves (\ref{eqHJB1simp})] we have that
%
%
\begin{eqnarray}\label{eqinterim3}\qquad
\phi_{\kappa}\lam(x)&\leq& \Ex_x^{\theta}\int_0^{t\wedge\hat
{\tau}_{\kappa}^n}e^{-\gamma s} L(\hat{X}^n(s),\hat{U}(s))\,ds+\Ex
_{x}^{\theta}e^{-\gamma(t\wedge\hat{\tau}_{\kappa}^n)}\phi
_{\kappa}\lam\bigl(\hat{X}^n(t\wedge\hat{\tau}_{\kappa}^n)\bigr)
\nonumber\\[-8pt]\\[-8pt]
&&{}-\Ex_{x}^{\theta}\sum_{i\in\I}\int_0^{t\wedge\hat{\tau
}_{\kappa}^n}
e^{-\gamma s} (\phi_{\kappa}\lam)_i(\hat{X}^n(s))\sigma_i^n(\hat
{X}^n(s),\hat{U}(s)) \,dB(s).\nonumber
\end{eqnarray}
Here, $\hat{\tau}_{\kappa}^n$ is as defined in Definition \ref
{definadmissiblesystemBrownian} and it is a stopping time with respect
to $(\mathcal{F}_t)$ because of the continuity of $\hat{X}^n$. We now
claim that
\[
\Ex_{x}^{\theta}\bigl[e^{-\gamma t\wedge\hat{\tau}_{\kappa
}^n}\phi_{\kappa}\lam\bigl(\hat{X}^n(t\wedge\hat{\tau}_{\kappa
}^n)\bigr)\bigr]\rightarrow0 \qquad\mbox{as } t\tinf.
\]
Indeed, as $\phi_{\kappa} \lam$ is bounded on $\mB_{\kappa}^n$, on
the event $\{\hat{\tau}_{\kappa}^n=\infty\}$ we have that
\[
e^{-\gamma(t\wedge\hat{\tau}_{\kappa}^n)}\phi_{\kappa}\lam
\bigl(X(t\wedge\hat{\tau}_{\kappa}^n)\bigr)\rightarrow0 \qquad\mbox{as }t\tinf.
\]
On the event $\{\hat{\tau}_{\kappa}^n<\infty\}$ we have $\hat
{X}^n(\hat{\tau}_{\kappa}^n)\in\partial\mB$ and, by the
definition of $\hat{\tau}_{\kappa}^n$, that $\phi_{\kappa}\lam
(\hX^n(\hat{\tau}_{\kappa}^n))=0$. The convergence in expectation
then follows from the bounded convergence theorem (using again the
boundedness of $\phi_{\kappa}\lam$ on $\mB_{\kappa}^n$). The last
term in (\ref{eqinterim3}) equals zero by the optional stopping
theorem.\vadjust{\goodbreak}

Letting $t\tinf$ in (\ref{eqinterim3}) and applying the monotone
convergence theorem, we then have
\[
\phi_{\kappa}\lam(x)\leq\Ex_x^{\theta} \biggl[\int_0^{\hat{\tau
}_{\kappa}^n}e^{-\gamma s} L(\hat{X}^n(s),\hat{U}(s))\,ds\biggr].
\]

Since the admissible system $\theta$ was arbitrary, we have that $\phi
_{\kappa}\lam(x)\leq\hat{V}\lam(x,\kappa)$. To show that this
inequality is actually an equality, let
%
%
\begin{equation}\label{eqratiofunchouce}
h_{\kappa
}^n(x)=e_{i^n(x)},
\end{equation}
where $e_{i^n(x)}$ is
as defined in the statement of the theorem.

The continuity of $\phi_{\kappa}^n$ guarantees that the function
$i^n(x)$ is Lebesgue measurable, and so is, in turn, $h_{\kappa
}^n(\cdot)$. Consider now the autonomous SDE:
%
%
\begin{equation}\label{eqautonomous}
\hX^n(t)=x+\int
_0^t \hat{b}^n(\hX^n(s))\,ds+\int_0^t \hat{\sigma}^n(\hX
^n(s))\,dB(s),
\end{equation}
where $\hat
{b}^n(y)=b^n(y,h_{\kappa}^n(y))$ and $\hat{\sigma}^n(y)=\sigma
^n(y,h_{\kappa}^n(y))$ on $\mB_{\kappa}^n$. Then, $\hat{b}^n$ and~%
$\hat{\sigma}^n$ are bounded and measurable on the bounded domain
$\mB_{\kappa}^n$. Also, as the matrix $\hat{\sigma}^n$ is diagonal
and the elements on the diagonal are strictly positive on $\mB_{\kappa
}^n$, it is positive definite there. Hence, a weak solution exists for
the autonomous SDE (see, e.g.,\vspace*{2pt} Theorem 6.1 of \cite
{krylov2008controlled}). In particular, there exists a probability\vadjust{\goodbreak}
space $(\tilde{\Omega},\mathcal{G},\tilde{\Pd})$, a filtration
$(\mathcal{G}_t)$ that satisfies the usual conditions, a Brownian
motion $B(t)$ and a continuous process $\hX^n$---both adapted to
$(\mathcal{G}_t)$, so that $\hX^n$ satisfies the autonomous SDE (\ref
{eqautonomous}). Finally, since~$\hX^n$ has continuous sample paths
and it is adapted, it is also progressively measurable (see,
e.g.,\vspace*{1pt}
Proposition 1.13 in \cite{KaS91}) and, by measurability of $h_{\kappa
}^n(\cdot)$, so is the process $\hat{U}(t)=h_{\kappa}^n(\hX^n(t))$.
Consequently, $\theta=(\tilde{\Omega},\mathcal{G},\mathcal
{G}_t,\tilde{\Pd},\hat{U},B)$ is an admissible system in the sense
of Definition \ref{definadmissiblesystemBrownian} and $\hX^n$ is the
corresponding controlled process.

To see that $\theta$ is optimal for the ADCP on $\mB_{\kappa}^n$,
note that for $s<\hat{\tau}_{\kappa}^n$, we have by the HJB equation
(\ref{eqHJB0}) that
\[
L(\hX^n(s),\hat{U}(s))+A_{\hat{U}(s)}\phi_{\kappa}\lam(\hX
^n(s))-\gamma\phi_{\kappa}\lam(\hX^n(s))=0.
\]
Applying It\^{o}'s rule as before, together with the bounded and
dominated convergence theorems, we then have that
\[
\phi_{\kappa}\lam(x)=\Ex_x^{\theta}\biggl[\int_0^{\hat{\tau
}_{\kappa}^n}e^{-\gamma s} L(\hX^n(s),\hat{U}(s))\,ds\biggr]
\]
and the proof is complete.

\section{The performance analysis of tracking policies}\label{sectracking}

This section shows that given an optimal Markov control policy for the
ADCP together with its associated tracking function $h_{\kappa
}^{*,n}$, the nonpreemptive tracking policy imitates, in a particular
sense, the performance of the Brownian system.
\begin{theorem} \label{thmSSC}
Fix $\kappa$ and $\kappa'<\kappa$ as well as a
sequence $\{(x^n,q^n),n\in\bbZ_+\}$ such that $(x^n,q^n)\in\mathcal
{X}^n$, and $|x^n-\nu n| \leq M \sqrt{n}$ for all $n$ and some
$M>0$. Let $\phi_{\kappa}^n$ and $h_{\kappa}^{*,n}$ be as in Theorem
\ref{thmBrownianverification} and define
\[
\psi^n(x,u)=L(x,u)+A^n_u \phi_{\kappa}^n(x)-\gamma\phi_{\kappa
}^n(x) \qquad\mbox{for }x\in\mB_{\kappa}^n,u\in\mathcal{U}.
\]
Let $U_h^n$ be the ratio control associated with the $h_{\kappa
}^{*,n}$-tracking policy and let $\mathbb{X}^n=(X^n,Q^n,Z^n,\check
{X}^n)$ be the associated queueing process with the initial conditions
$Q^n(0)=q^n$ and $\check{X}^n(0)=x^n-\nu n$ and define
\[
\tau_{\kappa',T}^n=\inf\{t\geq0\dvtx\check{X}^n(t)\notin\mB_{\kappa
'}^n\}\wedge T\log n.
\]
Then,
\[
\Ex\biggl[\int_0^{\tau_{\kappa',T}^n} e^{-\gamma s} |\psi
^n(\check{X}^n(s),U_h^n(s))-\psi^n(\check{X}^n(s),h_{\kappa
}^{*,n}(\check{X}^n(s)))|\,ds\biggr]\leq C\log^{k_0+3}n
\]
for a constant $C$ that does not depend on $n$.
\end{theorem}

Theorem \ref{thmSSC} is proved in the \hyperref[app]{Appendix}. The proof builds on the
gradient estimates in Theorem \ref{thmHJB1sol} and on a state-space
collapse-type result for certain sub-intervals of $[0,\tau_{\kappa',T}^n]$.
\begin{rem} \label{remSSC} Typically one establishes a stronger
state-space collapse result showing that the actual queue and the
desired queue values are close in supremum norm. The difficulty with
the former approach is that the tracking functions here are nonsmooth.
While it is plausible that one can smooth these functions appropriately
(as is done, e.g., in \cite{AMR02}), such smoothing might
compromise the optimality gap. Fortunately, the weaker integral
criterion implied by Theorem
\ref{thmSSC} suffices for our purposes.
\end{rem}

\section{Proof of the main result}\label{seccombining}

Fix $\kappa>0$ and let $\phi_{\kappa}^n$ be the solution to (\ref
{eqHJB1simp}) on $\mB_{\kappa}^n$ (see Theorem \ref{thmHJB1sol}).
We start with the following lemma where $b_i^n(\cdot,\cdot)$ and
$\sigma_i^n(\cdot,\cdot)$ are as in (\ref{eqnbmdefn30}) and
(\ref{eqnbmdefn3}), respectively.
\begin{lem} \label{lemito} Let $U^n$ be an admissible ratio control
and let $\mathbb{X}^n=(X^n,Q^n$, $Z^n,\check{X}^n)$ be the queueing
process associated with $U^n$. Fix $\kappa'<\kappa$ and $T>0$ and let
\[
\tau_{\kappa',T}^n=\inf\{t\geq0\dvtx\check{X}^n(t)\notin\mB_{\kappa
'}^n\}\wedge T\log n.
\]
Then, there exists a constant $C$ that does not depend on $n$ (but may
depend on $T$, $\kappa$ and $\kappa'$) such that
\begin{eqnarray*}
\Ex[e^{-\gamma\tau_{\kappa',T}^n}\phi_{\kappa}^n(\check
{X}^n(\tau_{\kappa',T}^n))]&\leq&\phi_{\kappa}^n(\check{X}^n(0))+
\Ex\biggl[\int_0^{\tau_{\kappa',T}^n} e^{-\gamma s}
A_{U^n(s)}^n\phi_{\kappa}^n(\check{X}^n(s))\,ds\biggr]\\
&&{}-\gamma\Ex
\biggl[\int_0^{\tau_{\kappa',T}^n} e^{-\gamma s} \phi_{\kappa
}^n(\check{X}^n(s))\,ds\biggr] +C\log^{k_1+1} n\\
& \leq&
\Ex[e^{-\gamma\tau_{\kappa',T}^n}\phi_{\kappa}^n(\check
{X}^n(\tau_{\kappa',T}^n))]+2C\log^{k_1+1}n.
\end{eqnarray*}
\end{lem}

We will also use the following lemma where $c=(c_1,\ldots, c_I)$ are
the cost coefficients (see Section \ref{secmodel}).
\begin{lem}\label{lemafterstop}
Let $(x^n,q^n)$ be as in the conditions of Theorem \ref
{thmmain}. Then, there exists a constant $C$ that does not depend on
$n$ such that
%
%
\begin{equation}\label{eqafterstop2}
\Ex_{x^n,q^n}^{U}\biggl[\int_{\tau_{\kappa
',T}^n}^{\infty} e^{-\gamma s} (e\cdot c)\bigl(e\cdot\check
{X}^n(s)\bigr)^+\,ds\biggr]\leq C\log^2n
\end{equation}
and
%
%
\begin{equation}\label{eqafterstop1}
\Ex_{x^n,q^n}^{U}[e^{-\gamma\tau_{\kappa
',T}^n}\phi
_{\kappa}^n(\check{X}^n(\tau_{\kappa',T}^n))]\leq C\log^2
n
\end{equation}
for all $n$ and any admissible ratio control
$U$.
\end{lem}

We postpone the proof of Lemma \ref{lemito} to the end of the section
and that of Lemma \ref{lemafterstop} to the \hyperref[app]{Appendix} and proceed now to
prove the main result of the paper.

\subsection*{\texorpdfstring{Proof of Theorem \protect\ref{thmmain}}{Proof of Theorem 2.1}} Let $h_{\kappa
}^{*,n}$ be the ratio function associated with the optimal Markov
control for the ADCP (as in Theorem \ref{thmHJB1sol}). Since $\kappa
$ is fixed we omit the subscript $\kappa$ and use $h^n=h_{\kappa
}^{*,n}$. Let $U_h^n$ be the ratio associated with the $h^n$-tracking policy.

The proof will proceed in three main steps. First, building on Theorem~\ref{thmSSC} we will show that
%
%
\begin{equation} \label{eqinterim2}
\Ex\biggl[\int_0^{\tau_{\kappa
',T}^n} e^{-\gamma s} L(\check{X}^n(s),U_h^n(s))\,ds\biggr] \leq\phi
_{\kappa}^n(\check{X}^n(0))+C\log^{k_0+3} n.
\end{equation}
Using Lemma \ref{lemafterstop}, this implies
%
%
\begin{eqnarray}\label{eqinterim13}
C\lam(x^n,q^n,U_h^n)&=&\Ex\biggl[\int_0^{\infty} e^{-\gamma s} L(\check
{X}^n(s),U_h^n(s))\,ds\biggr]\nonumber\\[-8pt]\\[-8pt]
&\leq&\phi_{\kappa}^n(\check
{X}^n(0))+C\log^{k_0+3} n.\nonumber
\end{eqnarray}
Finally, we will
show that for any ratio control $U^n$,
%
%
\begin{equation}\label{eqinterim303}
\phi_{\kappa}^n(\check
{X}^n(0))\leq\Ex\biggl[\int_0^{\infty} e^{-\gamma s} L(\check
{X}^n(s),U^n(s))\,ds\biggr]+C\log^{k_1+1}n,
\end{equation}
where we recall that $k_1=k_0-m$. In turn,
\[
V^n(x^n,q^n)\geq\phi_{\kappa}^n(x^n-\nu n)-C\log^{k_1+1} n \geq
C\lam(x^n,q^n,U_h^n)-2C\log^{k_1+1}n,
\]
which establishes the statement of the theorem.\vadjust{\goodbreak}

We now turn to prove each of (\ref{eqinterim2}) and (\ref{eqinterim303}).

\subsection*{\texorpdfstring{Proof of (\protect\ref{eqinterim2})}{Proof of (61)}} To simplify
notation we fix $\kappa>0$ throughout and let $h^n(\cdot
)=h_{\kappa}^{*,n}$. Using Lemma \ref{lemito} we have
%
%
\begin{eqnarray}\label{eqinterim1}
&&\Ex[e^{-\gamma\tau_{\kappa',T}^n}\phi_{\kappa}^n(\check
{X}^n(\tau_{\kappa',T}^n))]\nonumber
\\
&&\qquad\leq \phi_{\kappa}^n(\check{X}^n(0))+
\Ex\biggl[\int_0^{\tau_{\kappa',T}^n} e^{-\gamma s}
A_{U_h^n(s)}^n\phi_{\kappa}^n(\check{X}^n(s))\,ds\biggr]\\
&&\qquad\quad{}-\gamma\Ex\biggl[\int_0^{\tau_{\kappa',T}^n} e^{-\gamma s} \phi
_{\kappa}^n(\check{X}^n(s))\,ds\biggr]+ C\log^{k_1+1} n .\nonumber
\end{eqnarray}
From the definition of $h^n$ as a minimizer in the HJB equation we have that
\begin{eqnarray*}
0&=& \Ex\biggl[\int_0^{\tau_{\kappa',T}^n} e^{-\gamma s}
A_{h^n(\check{X}^n(s))}^n\phi_{\kappa}^n(\check{X}^n(s))\,ds\biggr]\\
&&{}-\gamma\Ex\biggl[\int_0^{\tau_{\kappa',T}^n} e^{-\gamma s} \phi
_{\kappa}^n(\check{X}^n(s))\,ds\biggr]\\
&&{}+\Ex\biggl[\int
_0^{\tau_{\kappa',T}^n} e^{-\gamma s} L(\check{X}^n(s),h^n(\check
{X}^n(s)))\,ds\biggr] .
\end{eqnarray*}
By Theorem \ref{thmSSC} we then have that
%
%
\begin{eqnarray}\label{eqinterim101}
C\log^{k_0+3}n&\geq& \Ex\biggl[\int_0^{\tau_{\kappa',T}^n}
e^{-\gamma s} A_{U_h^n(s)}^n\phi_{\kappa}^n(\check{X}^n(s))\,ds
\biggr]\nonumber\\[-2pt]
&&{} -\gamma\Ex\biggl[\int_0^{\tau_{\kappa',T}^n}
e^{-\gamma s}
\phi_{\kappa}^n(\check{X}^n(s))\,ds\biggr]\nonumber\\[-9pt]\\[-9pt]
&&{} +\Ex\biggl[\int_0^{\tau_{\kappa',T}^n} e^{-\gamma s} L(\check
{X}^n(s),U_h^n(s))\,ds\biggr]\nonumber\\[-2pt]
&\geq&0.\nonumber
\end{eqnarray}
Since $\phi_{\kappa}^n$ is nonnegative, combining (\ref{eqinterim1})
and (\ref{eqinterim101}) we have that
\[
\Ex\biggl[\int_0^{\tau_{\kappa',T}^n} e^{-\gamma s} L(\check
{X}^n(s),U_h^n(s))\,ds\biggr]\leq\phi_{\kappa}^n(\check
{X}^n(0))+C\log^{k_0+3} n,
\]
which concludes the proof of (\ref{eqinterim2}).

\subsection*{\texorpdfstring{Proof of (\protect\ref{eqinterim303})}{Proof of (63)}} We now show that
$V^n(x,q)\geq\phi_{\kappa}^n(\check{X}^n(0))-C\log^{k_1+1} n$. To
that end, fix an arbitrary ratio control $U^n$ and recall that by the
HJB equation,
\[
A_u^n\phi_{\kappa}^n(x)-\gamma\phi_{\kappa
}^n(x)+L(x,u)\geq0
\]
for all $u\in\mathcal{U}$ and $x\in\mB
_{\kappa}^n$. In turn, using the second inequality in Lemma \ref
{lemito} we have that
\begin{eqnarray*}
&&\Ex[e^{-\gamma\tau_{\kappa
',T}^n}\phi_{\kappa}^n(\check{X}^n(\tau_{\kappa',T}^n))
]\\[-2pt]
&&\qquad\geq \phi_{\kappa}^n(\check{X}^n(0))
-\Ex\biggl[\int_0^{\tau_{\kappa',T}^n} e^{-\gamma s} L(\check
{X}^n(s),U^n(s))\,ds\biggr]\\[-2pt]
&&\qquad\quad{} -2C\log^{k_1+1} n.
\end{eqnarray*}
Using Lemma \ref{lemafterstop}, we have, however, that
\[
\Ex[e^{-\gamma\tau_{\kappa',T}^n}\phi_{\kappa}^n(\check
{X}^n(\tau_{\kappa',T}^n))]\leq C\log^{2} n
\]
for a redefined constant $C$ so that
\begin{eqnarray*}C\log^{2} n &\geq& \phi
_{\kappa}^n(\check{X}^n(0))-
\Ex\biggl[\int_0^{\tau_{\kappa',T}^n} e^{-\gamma s} L(\check
{X}^n(s),U^n(s))\,ds\biggr]\\[-2pt]
&&{}-2C\log^{k_1+1} n\\&\geq& \phi_{\kappa
}^n(\check{X}^n(0))-
\Ex\biggl[\int_0^{\infty} e^{-\gamma s} L(\check
{X}^n(s),U^n(s))\,ds\biggr]\\[-2pt]
&&{}-2C\log^{k_1+1} n
\end{eqnarray*}
and, finally,
\[
\phi_{\kappa}^n(\check{X}^n(0))\leq\Ex\biggl[\int_0^{\infty}
e^{-\gamma s}
L(\check{X}^n(s),U^n(s))\,ds\biggr]+C\log^{k_1+1}n\vadjust{\goodbreak}
\]
for a redefined constant $C>0$. This concludes the proof of (\ref
{eqinterim303}) and of the theorem.

We end this section with the proof of Lemma \ref{lemito} in which the
following auxiliary lemma will be of use.
\begin{lem}\label{lemmartingales}
Fix $\kappa>0$ and an admissible ratio control $U^n$ and
let $\mathbb{X}\lam=(X\lam,Q\lam,Z\lam,\check{X}\lam)$ be the
corresponding queueing process. Let
\[
\tau_{\kappa,T}^n=\inf\{t\geq0\dvtx\check{X}^n(t)\notin\mB_{\kappa
}^n\}\wedge T\log n,
\]
and $(\check{W}_i^n,i\in\I)$ be as defined in (\ref{eqWtildedefin}).
Then, for each $i\in\I$, the process $\check
{W}_i^n(\cdot\wedge\tau_{\kappa,n}^n)$ is a square integrable
martingale w.r.t to the filtration $(\mathcal{F}_{t\wedge\tau
_{\kappa,T}^n}^n)$ as are the processes
\[
\mathcal{M}_i^n(\cdot)=\bigl(\check{W}_i^n(\cdot\wedge\tau_{\kappa
,T}^n)\bigr)^2-\int_0^{\cdot\wedge\tau_{\kappa,T}^n} (\sigma
_i^n(\check{X}^n(s),U^n(s)))^2\,ds
\]
and
\[
\mathcal{V}_i^n(\cdot)=\bigl(\check{W}_i^n(\cdot\wedge\tau_{\kappa
,T}^n)\bigr)^2-\sum_{s\leq\cdot\wedge\tau_{\kappa,T}^n} (\Delta\check
{W}_i^n(s))^2.
\]
\end{lem}

Lemma \ref{lemmartingales} follows from basic results on martingales
associated with time-changes of Poisson processes. The detailed proof
appears in the \hyperref[app]{Appendix}.

\subsection*{\texorpdfstring{Proof of Lemma \protect\ref{lemito}}{Proof of Lemma 6.1}}
Note that, as in
(\ref{eqcheckXdynamics}), $\check{X}^n$ satisfies
\[
\check{X}_i\lam(t)=\check{X}_i\lam(0)+\int_0^t b_i^n(\check
{X}^n(s),U^n(s))\,ds+\check{W}_i\lam(t),
\]
and is a semi martingale. Applying It\^{o}'s formula for
semimartingales (see, e.g., Theorem 5.92 in \cite{vandervaart}) we
have for all $t\leq\tau_{\kappa',T}^n$, that
\begin{eqnarray*}
e^{-\gamma t}\phi_{\kappa}^n(\check{X}^n(t))&=&\phi_{\kappa
}^n(\check{X}^n(0))\\
&&{}+
\sum_{s\leq t\dvtx |\Delta\check{X}^n(s)|> 0} e^{-\gamma s}[\phi
_{\kappa}^n(\check{X}^n(s))-\phi_{\kappa}^n(\check{X}^n(s-))]\\
&&{}-
\sum_{i\in\I}\sum_{s\leq t\dvtx |\Delta\check{X}^n(s)|> 0} e^{-\gamma
s}(\phi_{\kappa})_i^n(\check{X}^n(s))\Delta\check{X}_i^n(s)\\
&&{}+
\sum_{i\in\I}\int_0^t e^{-\gamma s} (\phi_{\kappa}^n)_i(\check
{X}^n(s-))b_i^n(\check{X}^n(s),U^n(s))\,ds\\
&&{}-\gamma\int_0^t
e^{-\gamma s} \phi_{\kappa}^n(\check{X}^n(s))\,ds
\end{eqnarray*}
and, after rearranging terms, that
\begin{eqnarray*}
&&e^{-\gamma t}\phi_{\kappa}^n(\check{X}^n(t))\\
&&\qquad=
\phi_{\kappa}^n(\check{X}^n(0))
+ \frac{1}{2}\sum_{i\in\I
}\sum_{s\leq t\dvtx |\Delta\check{X}^n(s)|>0}e^{-\gamma s} (\phi
_{\kappa}^n)_{ii}(\check{X}^n(s-))(\Delta\check{X}_i^n(s))^2\\
&&\qquad\quad{}+
\sum_{i\in\I}\int_0^t e^{-\gamma s} (\phi_{\kappa}^n)_i(\check
{X}^n(s-))b_i(\check{X}^n(s),U^n(s))\,ds\\
&&\qquad\quad{}+ C^n(t)-\gamma\int_0^t
e^{-\gamma s} \phi_{\kappa}^n(\check{X}^n(s))\,ds,
\end{eqnarray*}
where
\begin{eqnarray*} C^n(t)&=& \sum_{s\leq t\dvtx |\Delta\check
{X}^n(s)|>0}e^{-\gamma s}  \biggl[\phi_{\kappa}^n(\check
{X}^n(s))-\phi_{\kappa}^n(\check{X}^n(s-))\\
&&\hphantom{\sum_{s\leq t\dvtx |\Delta\check
{X}^n(s)|>0}e^{-\gamma s}  \biggl[}{}-\sum_{i\in\I}(\phi
_{\kappa}^n)_i(\check{X}^n(s-)) \Delta\check{X}_i^n(s) \\
&&\hphantom{\sum_{s\leq t\dvtx |\Delta\check
{X}^n(s)|>0}e^{-\gamma s}  \biggl[}{}-\frac
{1}{2}\sum_{i\in\I} (\phi_{\kappa}^n)_{ii}(\check
{X}^n(s-))(\Delta\check{X}_i^n(s))^2\biggr].
\end{eqnarray*}
Setting $t=\tau_{\kappa',T}^n$ as defined in the statement of the
lemma and taking expectations on both sides we have
%
%
\begin{eqnarray}\label{eqinterim404}
&&\Ex[e^{-\gamma\tau_{\kappa',T}^n}\phi_{\kappa}^n(\check
{X}^n(t))]\nonumber\\
&&\qquad=\phi_{\kappa}^n(\check{X}^n(0))
+
\sum_{i\in\I}\Ex\biggl[\int_0^{\tau_{\kappa',T}^n} e^{-\gamma s}
(\phi_{\kappa}^n)_i(\check{X}^n(s-)) b_i^n(\check
{X}^n(s),U^n(s))\,ds\biggr]\nonumber\\[-8pt]\\[-8pt]
&&\qquad\quad{}+\frac{1}{2}\sum
_{i\in\I}\Ex\biggl[\sum_{s\leq t\dvtx |\Delta\check
{X}^n(s)|>0}e^{-\gamma s} (\phi_{\kappa}^n)_{ii}(\check
{X}^n(s-))(\Delta\check{X}_i^n(s))^2\biggr]\nonumber\\
&&\qquad\quad{}+ \Ex
[C^n(\tau_{\kappa',T}^n)]-\gamma\Ex\biggl[\int_0^{\tau_{\kappa
',T}^n} e^{-\gamma s} \phi_{\kappa}^n(\check{X}^n(s))\,ds
\biggr].\nonumber
\end{eqnarray}

We will now examine each of the elements on the right-hand side of
(\ref{eqinterim404}). First, note that $\Delta\check
{X}_i^n(s)=\Delta\check{W}_i^n(s)$ and, in particular,
\begin{eqnarray*}
&&\Ex
\biggl[\sum_{s\leq\tau_{\kappa',T}^n\dvtx |\Delta\check
{X}^n(s)|>0}e^{-\gamma s} (\phi_{\kappa}^n)_{ii}(\check
{X}^n(s-))(\Delta\check{X}_i^n(s))^2\biggr]\\
&&\qquad=
\Ex
\biggl[\sum_{s\leq\tau_{\kappa',T}^n\dvtx |\Delta\check
{X}^n(s)|>0}e^{-\gamma s} (\phi_{\kappa}^n)_{ii}(\check
{X}^n(s-))(\Delta\check{W}_i^n(s))^2\biggr].
\end{eqnarray*}
Using the fact that $\mathcal{V}_i^n$, as defined in Lemma \ref
{lemmartingales}, is a martingale as well as the fact that $\phi
_{\kappa}^n(\check{X}^n(s))$ and its derivative processes are bounded
up to~$\tau_{\kappa'}^n$, we have that the processes
%
%
\begin{equation}\label{eqbarV}
\bar{\mathcal{V}}_i^n(\cdot):=\int_0^{\cdot\wedge
\tau
_{\kappa',T}^n} e^{-\gamma s}(\phi_{\kappa}^n)_{ii}(\check
{X}^n(s-))\,d\mathcal{V}_i^n(s)
\end{equation}
and
%
%
\begin{equation}\label{eqbarM}
\bar{\mathcal{M}}_i^n(\cdot):=\int_0^{\cdot\wedge
\tau
_{\kappa',T}^n} e^{-\gamma s}(\phi_{\kappa}^n)_{ii}(\check
{X}^n(s-))\,d\mathcal{M}_i^n(s)
\end{equation}
are themselves
martingales with $\bar{\mathcal{V}}_i^n(0)=\bar{\mathcal
{M}}_i^n(0)=0$ and in turn, by optional stopping, that
$\Ex[\bar{\mathcal{V}}_i^n(\tau_{\kappa',T}^n)]=\Ex[\bar
{\mathcal{M}}_i^n(\tau_{\kappa',T}^n)]$ (see,\vspace*{1pt} e.g., Lemma 5.45 in
\cite{vandervaart}). In turn, by the definition of $\mathcal
{M}_i^n(\cdot)$ and $\mathcal{V}_i^n(\cdot)$ we have
\begin{eqnarray*}
&&\Ex\biggl[\sum_{s\leq\tau_{\kappa',T}^n\dvtx
|\Delta\check{X}^n(s)|>0}e^{-\gamma s} (\phi_{\kappa
}^n)_{ii}(\check{X}^n(s-))(\Delta\check{W}_i^n(s))^2\biggr]\\
&&\qquad=\Ex\biggl[\int_0^t (\phi_{\kappa}^n)_{ii}(\check
{X}^n(s-))\,d(\check{W}_i^n(s))^2\biggr]\\
&&\qquad=\Ex\biggl[\int
_0^{\tau_{\kappa',T}^n} e^{-\gamma s} (\phi_{\kappa}^n)_{ii}(\check
{X}^n(s-))(\sigma_i^n(\check{X}^n(s),U^n(s)))^2\,ds\biggr].
\end{eqnarray*}
Plugging this back into (\ref{eqinterim404}) we have that
\begin{eqnarray*}
&&\Ex[e^{-\gamma\tau_{\kappa',T}^n}\phi_{\kappa}^n(\check
{X}^n(t))]\\
&&\qquad=\phi_{\kappa}^n(\check{X}^n(0))
+\sum_{i\in\I}\Ex\biggl[\int_0^{\tau_{\kappa',T}^n} e^{-\gamma s}
(\phi_{\kappa}^n)_{i}(\check{X}^n(s-)) b_i^n(\check
{X}^n(s),U^n(s))\,ds\biggr]\\
&&\qquad\quad{}+
\frac{1}{2}\sum_{i\in\I}
\Ex\biggl[\int_0^{\tau_{\kappa',T}^n} e^{-\gamma s} (\phi_{\kappa
}^n)_{ii}(\check{X}^n(s-))(\sigma_i^n(\check{X}^n(s),U^n(s)) )^2
\,ds\biggr]\\
&&\qquad\quad{}-\gamma\Ex\biggl[\int_0^{\tau_{\kappa',T}^n}
e^{-\gamma s} \phi_{\kappa}^n(\check{X}^n(s))\,ds\biggr]+ \Ex
[C^n(\tau_{\kappa',T}^n)],
\end{eqnarray*}
which, using the definition of $A_u^n$ in (\ref{eqgendefin}), yields
\begin{eqnarray*}
&&\Ex[e^{-\gamma\tau_{\kappa',T}^n}\phi_{\kappa}^n(\check
{X}^n(t))]\\
&&\qquad=\phi_{\kappa}^n(\check{X}^n(0))
+\Ex\biggl[\int_0^{\tau_{\kappa',T}^n} e^{-\gamma s}
A_{U^n(s)}^n\phi_{\kappa}^n(\check{X}^n(s))\,ds\biggr]\\
&&\qquad\quad{}-\gamma\Ex
\biggl[\int_0^{\tau_{\kappa',T}^n} e^{-\gamma s} \phi_{\kappa
}^n(\check{X}^n(s))\,ds\biggr]\\
&&\qquad\quad{}+ \Ex[C^n(\tau_{\kappa',T}^n)].
\end{eqnarray*}
To complete the proof it then remains only to show that there
exists a~constant $C$ such that
\[
|\Ex[C^n(\tau_{\kappa',T}^{n})]|\leq C\log
^{k_1+1} n.
\]
To that end, note that by Taylor's
expansion,
\begin{eqnarray*}\phi_{\kappa}^n(\check{X}^n(s))&=&\phi_{\kappa
}^n(\check{X}^n(s-))
+\sum_{i\in\I}(\phi_{\kappa}^n)_i(\check{X}^n(s-))\Delta
\check{X}_i^n(s)\\
&&{}
+\frac{1}{2}\sum_{i\in\I}(\phi_{\kappa}^n)_{ii}\bigl(\check
{X}^n(s-)+\eta_{\check{X}^n(s-)}\bigr)\Delta\check{X}_i^n(s),
\end{eqnarray*}
where
$\eta_{\check{X}^n(s-)}$ is such that $\check{X}^n(s-)+\eta_{\check
{X}^n(s-)}$ is between $\check{X}^n(s-)$ and $\check{X}^n(s-)+\Delta
\check{X}^n(s)$. In turn, adding and subtracting a term, we have that
%
%
\begin{eqnarray}\label{eqCninterim}
&&\phi_{\kappa}^n(\check{X}^n(s))-\phi_{\kappa
}^n(\check{X}^n(s-))-\sum_{i\in\I}(\phi_{\kappa}^n)_i(\check
{X}^n(s-)) \Delta\check{X}_i^n(s) \nonumber\\
&&\quad{}-\frac{1}{2}\sum_{i\in\I}(\phi_{\kappa}^n)_{ii}(\check
{X}^n(s-))(\Delta\check{X}_i^n(s))^2\\
&&\qquad = \sum_{i\in\I}\frac
{1}{2}\bigl((\phi_{\kappa}^n)_{ii}\bigl(\check{X}^n(s-)+\eta_{\check
{X}^n(s-)}\bigr) -(\phi_{\kappa}^n)_{ii}(\check{X}^n(s-))\bigr)(\Delta
\check{X}_i^n(s))^2.\nonumber
\end{eqnarray}
Since the jumps are of size $1$ and, with probability 1, there are no
simultaneous jumps, we have that $|\eta_{\check{X}^n(s-)}|\leq1$.
Adding the discounting, summing and taking expectations we have
%
%
\begin{eqnarray}\label{eqinterim2222}\quad
&&\Ex[C^n(t)]\nonumber\\
&&\qquad\leq\Ex\biggl[\sum_{s\leq t\dvtx
|\Delta\check{X}^n(s)|>0} e^{-\gamma s} \sum_{i\in\I}\frac
{1}{2}\max_{y\dvtx|y|\leq1}\bigl((\phi_{\kappa}^n)_{ii}\bigl(\check
{X}^n(s-)+y\bigr)\\
&&\qquad\quad\hspace*{163.6pt}{} -(\phi_{\kappa
}^n)_{ii}(\check{X}^n(s-))\bigr)(\Delta\check{X}_i^n(s))^2
\biggr],\nonumber
\end{eqnarray}
and a lower bound can be created by minimizing over $y$ instead of
maximizing. Using again the fact that $\Delta\check{X}_i^n(t)=\Delta
\check{W}_i^n(t)$ and that $\bar{\mathcal{M}}_i^n$ and $\bar
{\mathcal{V}}_i^n$ as defined in (\ref{eqbarM}) and (\ref{eqbarV})
are martingales, we have that
%
%
\begin{eqnarray}\label{eqinterim505}
\Ex[C^n(t)]&\leq&\Ex\biggl[\int_0^t \sum_{i\in\I
}\frac{1}{2}\max_{y\dvtx|y|\leq1}\bigl((\phi_{\kappa}^n)_{ii}\bigl(\check
{X}^n(s-)+y\bigr) \nonumber\\
&&\hspace*{92.3pt}{} -(\phi_{\kappa}^n)_{ii}(\check{X}^n(s-))\bigr)\\
&&\hspace*{60.6pt}{}\times(\sigma
_i^n(\check{X}^n(s),U^n(s)))^2\,ds\biggr].\nonumber
\end{eqnarray}
From (\ref{eqgenbound}) we have that
%
%
\begin{equation}\label{eqinterim1111}
\frac{1}{2}\biggl|\sum_{i\in\I}\bigl((\phi_{\kappa
}^n)_{ii}(y)-(\phi_{\kappa}^n)_{ii}(x)\bigr)(\sigma_i^n(x,u))^2
\biggr|\leq C\log^{k_1} n
\end{equation}
for all $u\in\mathcal{U}$ and $x,y\in\mB_{\kappa'}^n$ with
$|x-y|\leq1$. The proof is then concluded by plugging (\ref
{eqinterim1111}) into (\ref{eqinterim505}), setting $t=\tau_{\kappa
',T}^n$ and recalling that we can repeat all the above steps to obtain
a lower bound in (\ref{eqinterim505}) by replacing $\max_{y\dvtx|y|\leq
1}$ with $\min_{y\dvtx|y|\leq1}$ in~(\ref{eqinterim2222}).

\section{Concluding remarks}\label{secconclusions}

This paper proposes a novel approach for solving problems of dynamic
control of queueing systems in the Halfin--Whitt many-server
heavy-traffic regime. Its main contribution is the use of Brownian
approximations to construct controls that achieve optimality gaps that
are logarithmic in the system size. This should be contrasted with the
optimality gaps of size $o(\sqrt{n})$ that are common in the
literature on asymptotic optimality. A distinguishing feature of our
approach is the use of a \textit{sequence} of Brownian control problems
rather than a single (limit) problem. Having an entire sequence of
approximating problems allows us to perform a more refined analysis,
resulting in the improved optimality gap.

In further contrast with the earlier literature, in each of these
Brownian problems the diffusion coefficient depends on both the system
state and the control. Incorporating the impact of control on diffusion
coefficients allows us to track the performance of the policy better
but, at the same time, it leads to a more complex diffusion control
problem in which the associated HJB equation is fully nonlinear and
nonsmooth. For \textit{each} Brownian problem, we show that the HJB
equation has a sufficiently smooth solution that coincides with the
value function and that admits an optimal Markov policy. Most
importantly, we derive useful gradient estimates that apply to the
whole sequence and bound the growth rate of the gradients with the
system size. These bounds are crucial for controlling the approximation
errors when analyzing the original queueing system under the proposed
tracking control.

The motivating intuition behind our approximation scheme is that the
value functions of each queueing system and its corresponding Brownian
control problem ought to be close. In particular, the optimal control
for the Brownian problem should perform well for the queueing system.
Moreover, the optimal Markov control of the Brownian problem can be
approximated by a ratio (or tracking) control for the queueing system.
While these observations are ``correct'' at a high level, they need to
be qualified further. Our analysis underscores two sources of
approximation errors that need to be addressed in order to obtain the
refined optimality gaps.
First, the value function of the Brownian control problem may be
substantially different than that of the (preemptive) optimal control
problem for the queueing system. This difference must be quantified
relative to the system size, which we do indirectly through the
gradient estimates for the value function of the Brownian control
problem; this is manifested, for example, in the proof of Lemma \ref{lemito}.

The second source of error is in trying to imitate the optimal ratio
control of the approximating Brownian control by a tracking control in
the corresponding queueing system. The error arises because we insist
on having a~nonpreemptive control for the queueing system. Whereas
under a preemptive control, one may be able to rearrange the queues
instantaneously to match the tracking function of the Brownian system,
this is not possible with nonpreemptive controls. Instead, we carefully
construct and analyze the performance of the proposed nonpreemptive
tracking policy. In doing so, we prove that the tracking control
imitates closely the Brownian system with respect to a specific
integrated functional of the queueing dynamics (see Theorem \ref
{thmSSC} and Remark \ref{remSSC}). Here too, the gradient estimates
for the value function of the Brownian system play a key role.

While the focus of this paper has been a relatively simple model
to illustrate the key ideas behind our approach and the important steps
in the analysis, we expect that similar results can be established in
the cases of impatient customers, more general cost structures as well
as more general network structures.

As suggested by the preceding analysis, the viability of these
extensions and others will depend on whether it is possible to (a)
solve the
sequence of Brownian control problems and establish the necessary
gradient estimates and (b) establish the corresponding approximation
result for the nonpreemptive tracking control.

While we expect that the results of \cite{TandG} on fully nonlinear
elliptic PDEs can be invoked for the more general settings, extending
our analysis which builds on those results may not be always straightforward.
In particular, it is not immediately obvious how to generalize the
proof of the tracking result in Theorem \ref{thmSSC} to more general settings.

Nevertheless, we can make some observations about the extensions
mentioned above:
\begin{itemize}
\item\textit{General convex costs.} As discussed in the
\hyperref[sec1]{Introduction}, the analysis of the convex holding cost case will
probably be simpler as one tends to get ``interior'' solutions in that
case as opposed to the corner solutions in the linear cost case, which
causes nonsmoothness.
We expect that the enhanced smoothness (relative to the linear holding
cost case) will simplify the analysis of the HJB equations as well as
that of the tracking performance.

\item\textit{Abandonment.} Our starting point in the analysis is that,
among preemptive policies, work conserving policies are optimal. This
is not, in general, true when customers are impatient and may abandon
while waiting (see the discussion in Section 5.1 of \cite{AMR02}). As is
the case in \cite{AMR02}, our analysis will go through also for the
case of impatient customers provided that the cost structure is such
that work conservation is optimal among preemptive policies.
\item\textit{General networks.} Inspired by the generalization of
\cite{AMR02}, by Atar \cite{atar2005scheduling}, to tree-like networks,
we expect, for example, that such a generalization is viable in our
setting as well. Indeed, we expect that the analysis of the (sequence
of) HJB equations and the sequence of ADCPs be fairly similar for the
tree-like network setting. We expect that, in that more general
setting, it would be more complicated to bound the performance of the
tracking policies as in Theorem~\ref{thmSSC}.
\end{itemize}

\begin{appendix}\label{app}
\section*{Appendix}

\subsection*{\texorpdfstring{Proof of Lemma \protect\ref{lemexistenceofcontrolled}}{Proof of Lemma 3.1}}

Up to $\tau_{\kappa}^n$, both functions $b^n(\cdot,u)$ and $\sigma
^n(\cdot,u)$ are bounded and Lipschitz continuous (uniformly in $u$).
With these conditions satisfied, strong existence and uniqueness follow
as in Appendix D of \cite{FlemingSoner}. Specifically, strong
existence follows by successive approximations as in the proof of
Theorem~2.9 of \cite{KaS91} and uniqueness follows as in Theorem 2.5
there.

\subsection*{\texorpdfstring{Proof of Proposition \protect\ref{propsolPHJB}}{Proof of Proposition 4.1}} Fix
$\kappa>0$, $n\in\bbZ_+$ and $a>0$. Recall that (\ref{eqHJB2})
corresponds to finding $\phi_{\kappa,a}^n\in\mC^2(\mB)$ such that
%
%
\begin{equation}\label{eqHJB1min}
0=F_a[x,\phi_{\kappa,a}^n(x),D\phi_{\kappa
,a}^n(x),D^2\phi
_{\kappa,a}^n(x)],\qquad x\in\mB,
\end{equation}
and so that
$\phi_{\kappa,a}^n=0$ on $\partial\mB$ where \mbox{$F_a[\cdot]$} is as
defined in (\ref{eqFForm}). Then, Proposition~\ref{propsolPHJB}
will follow from Theorem 17.18 in \cite{TandG} upon verifying certain
conditions. The gradient estimates will also follow from \cite{TandG}
by carefully tracing some constants to identify their dependence on
$\kappa,n$ and $a$.

To that end, note that the function $F_a^i(x,z,p,r)$ [as defined in
(\ref{eqFdefin})] is linear in the $(z,p,r)$ arguments for all $k\in
\I$ and $x\in\mB$. In turn, this function is concave in these
arguments. Hence, to apply Theorem 17.18 of \cite{TandG} it remains to
establish that condition (17.53) of \cite{TandG} is satisfied for each
of these functions. In the following we suppress the constant $a>0$
from the notation. It suffices to show that there exist constants
$\underbar{\Lambda}\leq\bar{\Lambda}$ and $\eta$ such that
uniformly in $k\in\I$, $y=(x,z,p,r)\in\Gamma$, and $\xi\in\mathbb{R}^I$
%
%
\begin{eqnarray}\label{eqFcond1}
&\displaystyle
0< \underbar{\Lambda}|\xi^2|\leq\sum_{i,j}F^k_{i,j}[y] \xi_i\xi
_j\leq\bar{\Lambda} |\xi|^2,&\\
\label{eqFcond2}
&\displaystyle \max\{
|F^k_p[y]|,|F^k_z[y]|,|F^k_{rx}[y]|,|F^k_{px}[y]|,|F^k_{zx}[y]|\}
\leq\eta\underbar{\Lambda},&
\\
\label{eqFcond3}
&\displaystyle
\max\{ |F^k_x[y]|,|F^k_{xx}[y]|\}\leq\eta\underbar
{\Lambda}(1+|p|+|r|),&
\end{eqnarray}
where
\[
F^k_{i,j}(x,z,p,r)=\frac{\partial}{\partial r_{ij}}F^k(x,z,p,r),\qquad
F^k_{x_l}(x,z,p,r)=\frac{\partial}{\partial x_{l}}F^k(x,z,p,r)\vadjust{\goodbreak}
\]
and
\[
(F^k_{rx}(x,z,p,r))_{ilj}=\frac{\partial^2}{\partial r_{il}\,\partial
x_{j}}F^k(x,z,p,r).\vspace*{-2pt}
\]
The other cross-derivatives are defined similarly. We will show that we
can choose $\underbar{\Lambda}=\varepsilon_0 n$, $\bar{\Lambda}=
\varepsilon_1 n$, $\eta=\varepsilon_2$ for constants $\varepsilon
_0,\varepsilon_1$ and $\varepsilon_2$ that do not depend on $n$ and
$a$---this will be important in establishing the aforementioned
gradient estimates. To establish (\ref{eqFcond1}) note that, given
$\xi\in\mathbb{R}^I$,
%
%
\begin{equation}\label{eqFijDeriv}\quad
F^{k}_{ij}\xi_i\xi_j=\cases{
\frac{1}{2}\bigl(\lambda_i^n+\mu_i(\nu_in+x_i)\bigr)\xi
_i^2, &\quad for $i=j, i\neq k$,\vspace*{1pt}\cr
\frac{1}{2}\bigl(\lambda_i^n+\mu_i(\nu_in+x_i)\bigr)\xi_i^2 -\frac
{1}{2}f(e\cdot x), &\quad for $i=j=k$,\vspace*{1pt}\cr
0, &\quad otherwise.}\vspace*{-2pt}
\end{equation}
Hence,
\[
\sum_{i,j}F_{ij}^k \xi_i\xi_j=\frac{1}{2}\sum_{i\in\I}\bigl(\lambda
_i^n+\mu_i(\nu_in+x_i)\bigr)\xi_i^2-\frac{1}{2}f(e\cdot x)\xi_k^2.\vspace*{-2pt}
\]

Consequently, for $(x,z,r,p)\in\Gamma$ we have that
\[
\sum_{i,j}F_{ij}^k \xi_i\xi_j\leq I \bigl( \lambda+\mu_{\max
}n+\mu_{\max}\kappa\sqrt{n}\log^mn\bigr)\sum_{i\in\I}\xi_i^2+
\frac{1}{2}\kappa\sqrt{n}\log^m n\xi_k^2,\vspace*{-2pt}
\]
where $\mu_{\max}=\max_{k}\mu_k$. In particular, we can choose
$\varepsilon_1>0$ so that for all $n\in\bbZ$,
\[
\sum_{i,j}F_{ij}^k \xi_i\xi_j\leq\varepsilon_1 n.\vspace*{-2pt}
\]
To obtain the lower bound note that, for $y\in\Gamma$,
\[
\sum_{i,j}F_{ij}^k \xi_i\xi_j\geq\frac{1}{2}\Bigl(\min_{i\in\I
}\lambda_i^n+\min_{i\in\I}\mu_i\kappa\sqrt{n}\log^m n
\Bigr)\sum_{i\in\I}\xi_i^2
-\frac{1}{2}\xi_k^2\kappa\sqrt{n}\log^mn.\vspace*{-2pt}
\]
Hence, we can find $\varepsilon_0>0$ such that for all n,
\[
\sum_{i,j}F_{ij}^k \xi_i\xi_j\geq\varepsilon_0 n.\vspace*{-2pt}
\]
Note that above $\varepsilon_0$ and $\varepsilon_1$ can depend on
$\kappa$ but they do not depend on $n$ and $a$. Hence, we have
established (\ref{eqFcond1}) and we turn to (\ref{eqFcond2}). To
that end, note that
%
%
\begin{eqnarray}\label{eqFp}
F_{p_k}^k(x,z,p,r)&=&f(e\cdot x)+l_k\lam-\mu
_kx_k \quad\mbox{and }\nonumber\\[-9pt]\\[-9pt]
F_{p_i}^k(x,z,p,r)&=&l_i\lam-\mu_ix_i \qquad\mbox{for }
i\neq k.\nonumber\vspace*{-2pt}
\end{eqnarray}
Therefore,
\begin{eqnarray*} |F_{p}^k|&\leq& (e\cdot x)^+ +1+\sum_{i}|l_i\lam
+\mu_i x_i|\\[-2pt]
&\leq&
I\kappa\sqrt{n}\log^mn+1+I\max_{i}\bigl(|l_i\lam|+\mu_i\kappa\sqrt
{n}\log^m n\bigr),\vspace*{-2pt}\vadjust{\goodbreak}
\end{eqnarray*}
where we used the simple observation that $f(e\cdot x)\leq(e\cdot
x)^++1$. Clearly, we can choose $\varepsilon_2$ so that $|F_p^k|\leq
\varepsilon_2\varepsilon_0\sqrt{n}\log^m n$. Also $F_z^k=-\gamma$
and \mbox{$F_{zx}=0$} so that by re-choosing $\varepsilon_2$ large enough we
have $\max\{|F^k_z[y]|,|F^k_{zx}[y]|\}\leq\varepsilon_2\varepsilon
_0\sqrt{n}\log n$. Finally, by (\ref{eqFijDeriv}) we have that
\begin{eqnarray*}
F^k_{r_{ij}x_l}&=&0 \qquad\mbox{for } i\neq j,\\
F^k_{r_{ii}x_j}&=&0 \qquad\mbox{for } i\neq k, i\neq
j,\\
F^k_{r_{ii}x_i}&=&\frac{1}{2}\mu_i \qquad\mbox{for } i\neq k,\\
F^k_{r_{ii}x_i}&=&\frac{1}{2}\mu_i \qquad\mbox{for } i\neq
k,\\
F^k_{r_{kk}x_k}&=&\frac{1}{2}\mu_k-\frac{1}{2}\,\frac{\partial
}{\partial x_k}f(e\cdot x),\\
F^k_{r_{kk}x_j}&=&\frac{1}{2}\,\frac{\partial}{\partial x_k}f(e\cdot
x) \qquad\mbox{for } j\neq k.
\end{eqnarray*}
Thus,
\[
|F_{rx}^k|^2\leq\sum_{l\in\I}\frac{1}{2}\biggl|\frac{\partial
}{\partial x_l}f(e\cdot x)\biggr|^2+\frac{1}{2}\mu_{\max}\leq\frac
{1}{2}(1+\mu_{\max}),
\]
where we used the fact that $f(\cdot)$ is continuously differentiable
with Lipschitz constant $1$ (independently of $a$). Finally,
\[
F_{x_i}^k=\frac{\partial}{\partial x_i}f(e\cdot x)\biggl(c_k+\mu
_kp_k-\frac{1}{2}\mu_kr_{kk}\biggr)-\mu_ip_i+\frac{1}{2}\mu_ir_{ii},
\]
so that
%
%
\begin{equation}\label{eqFx}
|F_{x_i}^k|\leq|c_k|+\mu_k|p|+\tfrac{1}{2}\mu
_k|r|+\mu
_i|p|+\tfrac{1}{2}\mu_i |r|.
\end{equation}
Also, note that
\[
F_{x_ix_j}^k=\frac{\partial}{\partial x_i\,\partial x_j}f(e\cdot
x)\biggl(c_k+p_k-\frac{1}{2}r_{kk}\biggr),
\]
so that
\[
F_{x_ix_j}^k=\cases{
2\bigl[c_k+\mu_kp_k-\frac{1}{2}\mu_kr_{kk}\bigr], &\quad
if $|e\cdot x|\leq\frac{1}{4}$,\vspace*{2pt}\cr
0, &\quad otherwise.}
\]
Combining the above gives
\[
|F_{xx}^k|\leq\varepsilon_2\varepsilon_0(1+|p|+|r|)
\]
for suitably\vspace*{1pt} redefined $\varepsilon_2$ which concludes the proof that
the conditions (\ref{eqFcond1})--(\ref{eqFcond3}) hold with $\bar
{\Lambda}=\varepsilon_1n$, $\underbar{\Lambda}=\varepsilon_0n$
and $\eta=\varepsilon_2$. Having verified these conditions, the
existence and uniqueness of the solution $\phi_{\kappa,a}^n$ to
(\ref{eqHJB2}) now follows from Theorem~17.18 in \cite{TandG}.

To obtain the gradient estimates in (\ref{eqgradients})
we first outline how the solution~$\phi
_{k,a}^n$ is obtained in \cite{TandG} as a limit of solutions to
smoothed equations (we refer the reader to~\cite{TandG}, page 466, for
the more elaborate description). To that end, let
$F_a^i$ be as defined in (\ref{eqFdefin}) and for $y\in\Gamma$ define
%
%
\begin{equation}\label{eqGhPDE}
F^h[y]=G_{h}(F^1_a[y],\ldots,F^I_a[y]),
\end{equation}
where
\[
G_h(y)=h^{-I}\int_{\bar{y}\in\bbR^I}\rho\biggl(\frac{y-\bar
{y}}{h}\biggr)G_0(\bar{y})\,d\bar{y}
\]
and $G_0(x)=\min_{i\in\I}x_i$ and $\rho(\cdot)$ is a mollifier on
$\bbR^I$ (see \cite{TandG}, page 466). $F^h$ satisfies all the bounds
in (\ref{eqFcond1})--(\ref{eqFcond3}) uniformly in $h$; cf. \cite
{TandG}, page 466. Then, there exists a unique solution $u^h$ for the equations
%
%
\begin{equation} \label{eqGhPDE2}
F^h[u^h]=0
\end{equation}
on $\mB_{\kappa}^n$ with
$u^h=0$ on $\partial\mB_{\kappa}^n$.

The solution $\phi_{\kappa,a}^n$ is now obtained as a limit of $\{
u^h\}$ in the space $C_*^{2,\alpha}(\mB)$ as defined in (\ref
{eqCstar}). Moreover, since the gradient bounds are shown in \cite
{TandG} to be independent of $h$, it suffices for our purposes to fix
$h$ and focus on the construction of the gradient bounds.

Our starting point is the bound at the bottom of page 461 of \cite
{TandG} by which
%
%
\begin{equation}\label{eqonemoreinterim}
|u^h|^*_{2,\alpha,\mB_{\kappa}^n}\leq\check
{C}(a,n)(1+|u^h|^*_{2,\mB_{\kappa}^n}),
\end{equation}
where\vspace*{-1pt} $|u^h|^*_{2,\mB_{\kappa}^n}=\sum_{j=0}^2 [u^h]_{j,\mB}^*$ and
$[\cdot]_{j,\mB}^*, j=0,1,2$, are as defined in Section~\ref{secADCP}.
The constant $\alpha(a,n)$ depends only on the number of
classes $I$ and on $\bar{\Lambda}/\underbar{\Lambda}$ (see~\cite
{TandG}, top of page 461) and this fraction equals, in our context, to
$\varepsilon_1/\varepsilon_0$ and is thus constant and independent of
$n$ and $a$.

We will\vspace*{1pt} address the constant $\check{C}(a,n)$ shortly. We first argue
how one proceeds from (\ref{eqonemoreinterim}). Fix $0<\delta< 1$,
let $\epsilon=\delta/\check{C}(a,n)$ and $C(\epsilon)=2/(\epsilon
/8)^{1/\alpha}$ (see \cite{TandG}, top of page 132). Then, applying
an interpolation inequality (see~\cite{TandG}, bottom of page 461 and
Lemma 6.32 on page 130), it is obtained that
\[
|u^h|^*_{2,0,\mB_{\kappa}^n}\leq C(\epsilon)|u^h|^*_{0,\Omega
}+\epsilon|u^h|^*_{2,\alpha,\mB_{\kappa}^n}.
\]
Plugging this back into (\ref{eqonemoreinterim}) one then has
\[
|u^h|^*_{2,\alpha,\mB_{\kappa}^n}\leq\check{C}(a,n)\biggl(1+\bar
{C}\check{C}(a,n)^{1/\alpha}|u^h|^*_{0,\mB_{\kappa}^n}
+\frac{\delta}{\check{C}(a,n)}|u^h|^*_{2,\alpha,\mB_{\kappa
}^n}\biggr)
\]
for a constant $\bar{C}$ that depends only on $\delta$ and $\alpha$.
In turn,
\[
|u^h|^*_{2,\alpha,\mB_{\kappa}^n}\leq\bar{C}(\check
{C}(a,n))^{1+1/\alpha}|u^h|^*_{0,\mB_{\kappa}^n}
\]
for a constant $\bar{C}$ that\vadjust{\goodbreak} does not depend on $a$ or $n$.

Hence, to obtain the required bound in (\ref{eqgradients})
it remains only to\break bound~$\check{C}(a,n)$. Following \cite{TandG},
building on equation (17.51) of \cite{TandG}, $\check{C}(a,n)$ is the
(minimal) constant
that satisfies
%
%
\begin{equation}\label{eqCandefin}
C(1+M_2)(1+\tilde{\mu}R_0+\bar{\mu}R_0^2)\leq
\check{C}(a,n)(1+|u^h|^*_{2,\mB}),
\end{equation}
where (as stated in \cite{TandG}, bottom of page 460) the (redefined)
constant $C$ depends only on the number of class $I$ and on $\bar
{\Lambda}/\underbar{\Lambda}=\varepsilon_1/\varepsilon_0$.
The constants~$\tilde{\mu}$ and~$\bar{\mu}$ are defined in \cite{TandG}
and we will explicitly define them shortly. Here one should not
confuse $\bar{\mu}$ with the average service rate in our system. In
what follows $\bar{\mu}$ will only be used as the constant in \cite
{TandG}. We now bound constants~$\tilde{\mu}$ and~$\bar{\mu}$.
These are defined by
\begin{eqnarray*}
\tilde{\mu}&=&\frac{D_0}{\underbar{\Lambda
}(1+M_2)},\qquad
\bar{\mu} =\frac{C(I)}{\underbar{\Lambda}}
\biggl(\frac{A_0^2}{\underbar{\Lambda}\epsilon}+\frac{B_0}{1+M_2}\biggr),
\\
D_0&=&\sup_{x,y\in\mB}\{|F^h_x(y,u^h(y),Du^h(y),D^2u^h(x))|
\\
&&\hphantom{\sup_{x,y\in\mB}\{}{}+|F^h_z(y,u^h(y),Du^h(y),D^2u^h(x))||Du^h(y)|\\
&&\hphantom{\sup_{x,y\in\mB}\{}{}+|F^h_p(y,u^h(y),Du^h(y),D^2u^h(x))||D^2u^h(y)|\},
\\
A_0&=&\sup_{\mB}\{|F^h_{rx}|+|F^h_p|\}, \\
B_0&=& \sup_{\mB} \{|F_{px}||D^2 u^h|+ |F_z||D^2
u^h|+|F_{zx}||Du^h|+|F_{xx}|\},
\end{eqnarray*}
where $C(I)$ is a constant that depends only on the number of classes
$I$, $\epsilon\in(0,1)$ is arbitrary and fixed (independent of $n$
and $a$) and $M_2=\sup_{\mB}|D^2u^h|$. The constants $\bar{\mu}$,
$\tilde{\mu}$ and $M_2$ are defined in \cite{TandG}, pages 456--460,
and $A_0$ and~$B_0$ are as on page 461 there.

We note that $F^h_z$ is a constant, $F^h_p$ is bounded by $\bar
{C}\sqrt{n}\log^m n$ for some constant $\bar{C}$ [see (\ref
{eqgradients})] that depends only on $\kappa$ and, by (\ref{eqFx}),
$|F_x^h|\leq\varepsilon_2\varepsilon_0(1+|p|+|r|)$. In turn,
$D_0\leq4\varepsilon_2\varepsilon_0\sqrt{n}\log^m n\sup_{\mB
}(1+|Du^h|+|D^2u^h|)$. Arguing similarly for $A_0$ and $B_0$ we find
that there exists a constant $\bar{C}$ (that does not depend on $n$
and $a$) such that
\[
A_0\leq\bar{C}\sqrt{n}\log^mn\quad\mbox{and}\quad B^0\leq\bar{C} \sup
_{\mB}(1+|Du^h|+|D^2u^h|),
\]
which in turn implies the existence of a redefined constant $\bar{C}$
such that
\[
\tilde{\mu}\leq\frac{\bar{C}\log^{m} n}{\sqrt{n}(1+M_2)}\sup
_{\mB}(1+|Du^h|+|D^2u^h|)
\]
and
\[
\bar{\mu}\leq\frac{\bar{C}\log^{2 m} n}{n}+\frac{\bar
{C}}{n(1+M_2)}\sup_{\mB}(1+|Du^h|+|D^2u^h|).
\]
The proof of the bound is concluded by plugging these back into (\ref
{eqCandefin}) and setting $R_0=\kappa\sqrt{n}\log^m n$ there to get
that
\[
\check{C}(a,n)\leq C\log^{4m(1+{1}/{\alpha})}n
\]
for some $C$ that does not depend on $a$ and $n$.

The constant $\tilde{C}$ on the right-hand side of (\ref{eqgradients})
(which can depend on $n$ but does not depend on $a$) is
argued as in the proof of Theorem 17.17 in~\cite{TandG} and we
conclude the proof by noting that the global Lipschitz constant (that
we allow to depend on~$n$) follows from Theorem 7.2 in
\cite{trudinger1983fully}.

We next turn to proof of Theorem \ref{thmSSC}. First, we will
explicitly construct the queueing process under the $h$-tracking policy
and state a lemma that will be of use in the proof of the theorem.
Define $A_i^n(t)=\mN_i^a(\lambda_i^n t)$ so that~$A_i^n$ is the
arrival process of class-$i$ customers. Given a ratio control $U^n$ and
the associated queueing process $\mathbb{X}^n=(X^n,Q^n,Z^n,\check{X}^n)$,
$\check{W}^n$ is as defined in~(\ref{eqWtildedefin}). Also, we define
\[
D^n(t)=\sum_{i\in\I}\mN_i^d\biggl( \mu_i\int_0^t Z_i\lam(s)\,ds
\biggr).
\]
That is, $D^n(t)$ is the total number of service completions by time
$t$ in the $n$th system.

For the construction of the queueing process under the tracking policy
we define a family of processes
$\{\mathcal{A}_{i,\mH}^n, i\in\I,\mathcal{H}\subset\I\}$ as follows:
let $\{\xi_{\mK}^l; l\in\bbZ_+,\mK\subset\I\}$ be a family of i.i.d
uniform $[0,1]$ random variables independent of
$\bar{\mathcal{F}}_{\infty}$ as defined in~(\ref{eqcheckFdefin}). For
each $\mathcal{K}\subset\I$, define the processes
$(\mathcal{A}_{i,\mH}^n, i\in\I)$ by
%
%
\begin{equation}\label{eqmAdefin}
\mathcal{A}_{i,\mH}^n(t)=\sum_{l=1}^{D^n(t)}1\biggl\{\frac{\sum_{k<
i,k\in\mH}\lambda_k}{1\vee\sum_{k\in\mH}\lambda_k} < \xi_{\mK
}^l\leq
\frac{\sum_{k\leq i,k\in\mH}\lambda_k}{1\vee
\sum_{k\in\mH}\lambda_k}\biggr\}.
\end{equation}

We note that for any strict subset $\mH\subset\I$ and $i\in\mK$,
the probability that a jump of $D^n(t)$ results in a jump of $\mathcal
{A}_{i,\mH}^n$ is equal to $\lambda_i^n/\sum_{k\in\mH}\lambda
_k^n=a_i/\sum_{k\in\mH}a_k$ and is strictly greater than $\lambda
_i^n/\sum_{k\in\I}\lambda_k^n=a_i$. We define
%
%
\begin{equation}\label{eqepsilondefin}
\epsilon_i=\min
_{\mathcal{H}\subset\I}a_i-\frac{a_i}{\sum_{k\in\mH}a_k},
\end{equation}
and note that $\epsilon_i>0$ by our assumption
that $a_i>0$ for all $i\in\I$ (see Section~\ref{secmodel}). Let
$\bar{\epsilon}=\min_{i}\epsilon_i/4$.

Note that at time intervals in which $i\in\mK(\cdot)=\mH$ (see
Definition \ref{defintracking}) for some $\varnothing\neq\mH\subset
\I$, the process $\mathcal{A}_{i,\mH}^n$ jumps with probability
$\lambda_i^n/\sum_{k\in\mH}\lambda_k$ whenever a server becomes
available (i.e., upon a jump of $D^n$). In turn, we will use the
processes $\{\mathcal{A}_{i,\mH}^n, i\in\I,\mathcal{H}\subset\I\}
$ to generate (randomized) admissions to service of class-$i$ customers
under the $h$-tracking policy.

More specifically, under the $h$-tracking policy (see Definition \ref
{defintracking}) a customer from the class-$i$ queue enters service in
the following events:
\begin{longlist}
\item A class-$i$ customer that arrives at time $t$ enters
service immediately if there are idle servers, that is, if $(e\cdot
\check{X}\lam(t-))^{-}>0$.
\item If a server becomes available at time $t$ (corresponding
to a jump of~$D^n$) and $t$ is such that $i\in\mathcal{K}(t-)=\mH
\subset\I$, then a customer from the class-$i$ queue is admitted to
service at time $t$ with probability $\lambda_i^n/\sum_{k\in\mH
}\lambda_k$. This admission to service corresponds to a jump of the
process $\mathcal{A}_{i,\mH}^n$ as defined in (\ref{eqmAdefin}).
\item If a server becomes available at time $t$ (corresponding
to a jump of~$D^n$) and $t$ is such that $\mathcal{K}(t-)=\varnothing$
and $i=\min\{k\in\I\dvtx Q_i^n(t)>0\}$, then a~class-$i$ customer is
admitted to service.
\end{longlist}
Formally, the queueing process $\mathbb{X}\lam=(X\lam,Q\lam,Z\lam
,\check{X}\lam)$ satisfies
\begin{eqnarray*} Z_i\lam(t)&=&Z_i\lam(0)+\int_0^t 1\bigl\{\bigl(e\cdot
\check{X}^n(s)\bigr)^->0\bigr\}\,dA_i^n(s)\\
&&{}+\sum_{\mH\subset\I}\int
_0^t1\{i\in\mathcal{K}(s-),\mathcal{K}(s-)=\mathcal{H}\}\,
d\mathcal{A}_{i,\mH}\lam(s)\\
&&{}+ \int_0^t 1\bigl\{\mK
(s-)=\varnothing,i=\min\{k\in\I\dvtx Q_k^n(s-)>0\}\bigr\}\,dD^n(s)\\
&&{}-\mN
_i^d\biggl( \mu_i\int_0^t Z_i\lam(s)\,ds \biggr),\qquad i\in\I,\\
X_i\lam(t)&=&X_i\lam(0)+A_i^n(t)-\mN_i^d\biggl(\mu_i \int_0^t
Z_i\lam(s) \,ds \biggr),\qquad i\in\I, \\ Q_i\lam(t)&=&X_i\lam(t)-Z_i\lam
(t),\qquad i\in\I.
\end{eqnarray*}

The second, third and fourth terms on the right-hand side of the
equation for $Z_i^n$ correspond, respectively, to the events described
by items (i)--(iii) above. Finally, $\check{X}\lam$ is defined from
$X\lam$ as in (\ref{eqtildeXdefin}). The fact that the above system
of equations has a unique solution is proved by induction on arrival
and service completions times (see, e.g., the proof of Theorem 9.2 of
\cite{MMR98}). Clearly, $\mathbb{X}^n$ satisfies (\ref
{eqdynamics2})--(\ref{eqnon-negativity2}) with $U_i\lam$ there
constructed from $Q\lam$ using~(\ref{eqUQmap}).

We note that, with this construction, the tracking policy is admissible
in the sense of Definition \ref{definadmissiblecontrols}. Also, it
will be useful for the proof of Theorem~\ref{thmSSC} to note that
with this construction, if $[s,t]$ is an interval such that \mbox{$i\in
\mathcal{K}(u)\subset\I$} for all $u\in[s,t]$ then
%
%
\begin{equation}\label{eqQtrack}\qquad
Q_i^n(t)-Q_i^n(s)=A_i^n(t)-A_i^n(s)-\sum_{\mH\subset
\I}\int
_s^t1\{\mathcal{K}(u-)=\mathcal{H}\}\,d\mathcal{A}_{i,\mH
}\lam(u).
\end{equation}

Before proceeding to the proof of Theorem \ref{thmSSC} the following
lemma provides preliminary bounds for arbitrary ratio controls.
\begin{lem}\label{lemstrongappbounds}
Fix $\kappa,T>0$ and a ratio control $U^n$, let $\mathbb
{X}^n=(X^n,Q^n,\allowbreak Z^n,\check{X}^n)$ be the associated queueing process
and define
\[
\tau_{\kappa,T}^n=\inf\{t\geq0\dvtx\check{X}^n(t)\notin\mB_{\kappa
}^n\}\wedge T\log n.
\]
Then, there exist constants $C_1,C_2,K_0>0$ (that depend on $T$ and
$\kappa$ but that do not depend on $n$ or on the ratio control $U^n$)
such that for all $K>K_0$ and all $n$ large enough,
%
%
\begin{eqnarray}\label{eqstrapp1}
&&\Pd\Bigl\{\sup_{0\leq t\leq2T\log n}|\check
{W}^n(t)|> K\sqrt
{n}\log n\Bigr\}\nonumber\\[-8pt]\\[-8pt]
&&\qquad\leq C_1e^{-C_2K\log n},\nonumber\\
\label{eqstrapp2}
&&\Pd\bigl\{|\check{X}^n(t)-\check{X}^n(s)|> \bigl((t-s)+(t-s)^2
\bigr)K\sqrt{n}\log n+K\log n,
\nonumber\\
&&\hspace*{163pt}
\mbox{ for some } s<t\leq2T\log n\bigr\}\\
&&\qquad\leq C_1 e^{-C_2K\log n},\nonumber
\end{eqnarray}
%
%
\begin{eqnarray}\label{eqstrapp3}
&&\Pd\{|A_i^n(t)-A_i^n(s)-\lambda_i^n(t-s)|
> \bar{\epsilon}n(t-s)+K\log n\nonumber\\
&&\hspace*{139.3pt} \mbox{ for some }
s<t\leq\tau_{\kappa,T}^n\}\\
&&\qquad\leq C_1e^{-C_2K\log n},\qquad i\in\I,\nonumber\\
\label{eqstrapp4}
&&\Pd\biggl\{
\biggl|D^n(t)-D^n(s)-\sum_{i}\mu_i\nu_i n(t-s)\biggr|
> \bar{\epsilon}n(t-s)+K\log n\nonumber\\
&&\hspace*{174pt}\mbox{ for some } s<t\leq\tau_{\kappa,T}^n\biggr\}
\\
&&\qquad\leq C_1e^{-C_2K\log n},\nonumber\\
\label{eqstrapp5}
&&\Pd\biggl\{\mathcal{A}_{i,\mH}^n(t)-\mathcal{A}_{i,\mH
}^n(s)-\lambda_i^n(t-s)\leq\frac{\epsilon_i}{2} n(t-s)-K\log
n \nonumber\\
&&\hspace*{156.6pt}\mbox{ for some } s<t\leq\tau
_{\kappa,T}^n\biggr\}\\
&&\qquad\leq C_1e^{-C_2K\log n},\qquad i\in\I,\mK\subset\I
.\nonumber
\end{eqnarray}
\end{lem}
\begin{pf}
Equation (\ref{eqstrapp1}) follows from strong approximations (see,
e.g., Lem\-ma~2.2. in \cite{csorgo}) and known bounds on the supremum of
Brownian motion (see, e.g., equation 2.1.53 in \cite{csorgo}).
Equation (\ref{eqstrapp2}) then follows using this bound together
with (52) in \cite{AMR02} but with $X^n(t)-X^n(s)$ instead of
$X^n(t)$ (in the notation of \cite{AMR02} $\check{W}^n$ is $\hat
{W}^n$). Equations (\ref{eqstrapp3})--(\ref{eqstrapp5}) follow by
carefully constructing and bounding the increments. We outline the
proof of (\ref{eqstrapp3}) and the others follow similarly. To
that\vadjust{\goodbreak}
end, note that given $K$ and for all $n$ large enough
\begin{eqnarray*}
&& \{|A_i^n(t)-A_i^n(s)-\lambda_i^n(t-s)|\leq
\bar
{\epsilon}n(t-s)+K\log n, \\
&&\hspace*{109pt}
\mbox{ for all }
0\leq s\leq t\leq2T\log n\} \\
&&\qquad \supseteq
\biggl\{\max_{l\leq N_i^n}\max_{j\geq0\dvtx j\/\log^l n\leq3T\log
n}\frac{|A_i^{j,l,n}-\lambda_i^n/\log^l n|}{\sqrt{\lambda_i^n/\log
^l n}}\leq K\sqrt{\log n}\biggr\},
\end{eqnarray*}
where
$A_i^{j,l,n}=A_i^n((j+1)/\log^ln)-A_i^n(j/\log^ln)$ and $N=\max\{l\dvtx
\log^ln \leq\lambda_i^n/\log n\}$. Indeed, given an interval $[s,t)$
we can construct it from smaller intervals. Starting with $l=0$, we fit
as many intervals of size $1$ into $[s,t)$, we then continue to fit as
many intervals of size $1/\log n$ to the uncovered part of the interval
and continue sequentially in $l$. We omit the simple and detailed
construction. Note that with such construction, given an interval
$[s,t)$, its covering uses at most $\log n$ intervals of size $\log^l
n$ for each $l\geq0$. Also, note that $N_i^n\leq C\log n$ for all $n$
and some constant $C$. From here, using strong approximations (or
bounds for Poisson random variables as in \cite{glynn1987upper}) we
have, for each $j$ and $l$, that
\[
\Pd\biggl\{\frac{|A_i^{j,l,n}-\lambda_i^n/\log^l n|}{\sqrt{\lambda
_i^n/\log^l n}}> K\sqrt{\log n}\biggr\}\leq C_1e^{-C_2 K\log n}.
\]
Since the number of intervals considered is of the order of $n\log n$,
the bound follows with redefined constants $C_1$ and $C_2$.
\end{pf}

\subsection*{\texorpdfstring{Proof of Theorem \protect\ref{thmSSC}}{Proof of Theorem 5.1}} Since $\kappa$
is fixed throughout we use $h^n(\cdot)=h_{\kappa}^{*,n}(\cdot)$. As
in the statement of the theorem, let
\[
\psi^n(x,u)=L(x,u)+A^n_u\phi_{\kappa}^n(x)-\gamma\phi_{\kappa
}^n(x)\qquad \mbox{for } x\in\mB_{\kappa}^n, u\in\mathcal{U}
,\vadjust{\goodbreak}
\]
so that by the definition of $A^n_u(x)$ we have
%
%
\begin{eqnarray}\qquad
\psi^n(x,u)&=&-\gamma\phi_{\kappa}^n(x)\nonumber\\
&&{} + (e\cdot
x)^+\cdot\sum_{i\in\I}u_i\biggl\{ c_i+\mu_i(\phi_{\kappa
}^n)_i(x)-\mu_i\frac{1}{2}(\phi_{\kappa}^n)_{ii}(x)\biggr\}
\\
&&{}+\sum_{i\in\I} (l_i\lam-\mu_ix_i)(\phi_{\kappa}^n)_i(x)
+\frac{1}{2}\sum_{i\in\I} \bigl(\lambda_i^n+\mu_i(\nu
_in+x_i)\bigr)(\phi_{\kappa}^n){ii}(x).\nonumber
\end{eqnarray}
Defining, as before,
\[
M_{i}^n(z)=c_{i}+\mu_i(\phi_{\kappa}^n)_{i}(z)-\tfrac{1}{2}\mu
_i(\phi_{\kappa}^n)_{ii}(z),
\]
we have that
\[
\psi^n(x,u)-\psi^n(x,v)=(e\cdot x)^+\biggl(\sum
_{i\in\I}v_iM_i^n(x)-\sum_{i\in\I}u_iM_i^n(x)\biggr).
\]
Let $U^n$ be the ratio control associated with the $h^n$-tracking
policy, let $\mathbb{X}^n=(X^n,Q^n,Z^n,\check{X}^n)$ be the
associated queueing process and define
%
%
\begin{eqnarray}\label{eqpsicheckdefin}
\check{\psi}^n(s)&=&\psi^n(\check
{X}^n(s),U^n(s))-\psi^n(\check{X}^n(s),h^n(\check{X}^n(s)))\nonumber
\\
&=&
\bigl(e\cdot\check{X}^n(s)\bigr)^+\sum_{i\in\I}h_i^n(\check
{X}^n(s))M_i^n(\check{X}^n(s)) \\
&&{}-\bigl(e\cdot
\check{X}^n(s)\bigr)^+ \sum_{i\in\I}U_i^n(s)M_i^n(\check
{X}^n(s)).\nonumber
\end{eqnarray}
Recall that, by construction, $Q_i^n(s)=(e\cdot\check
{X}^n(s))^+U_i^n(s)$ so that (\ref{eqpsicheckdefin}) can be
re-written as
\begin{eqnarray*} \check{\psi}^n(s)&=&\psi^n(\check
{X}^n(s),U^n(s))-\psi^n(\check{X}^n(s),h^n(\check{X}^n(s)))
\\&=&
\bigl(e\cdot\check{X}^n(s)\bigr)^+\sum_{i\in\I}h_i^n(\check
{X}^n(s))M_i^n(\check{X}^n(s)) \\
&&{}-\sum_{i\in\I
}Q_i^n(s)M_i^n(\check{X}^n(s)).
\end{eqnarray*}
The theorem will be proved if we show that
%
%
\begin{equation}\label{eqwhatneed}
\Ex\biggl[\int_0^{\tau_{\kappa',T}^n}e^{-\gamma s}
|\check
{\psi}^n(s)|\,ds\biggr]\leq C\log^{k_0+3}n.
\end{equation}
To
that end, define a sequence of times $\{\tau_{l}^n\}$ as follows:
\[
\tau_{l+1}^n=\inf\{t> \tau_l^n\dvtx h^n(\check{X}^n(t))\neq h^n(\check
{X}^n(\tau_l^n))\}\wedge\tau_{\kappa',T}^n \qquad\mbox{for } l\geq0,
\]
where $\tau_0^n=\eta^n\wedge\tau_{\kappa',T}^n$ and
%
%
\begin{equation}\label{eqetandefin}
\eta
^n=t_0\frac{\log^mn}{\sqrt{n}}
\end{equation}
for\vspace*{1pt}
$t_0=4\kappa/\epsilon_i$ with $
\epsilon_i=\min_{\mathcal{H}\subset\I}a_i-\frac{a_i}{\sum_{k\in
\mH}a_k}$ as in (\ref{eqepsilondefin}). Finally, we define $r^n=\sup
\{l\in\bbZ_+\dvtx \tau_{l}^n\leq\tau_{\kappa',T}^n\}$ and set $\tau
_{r^n+1}^n=\tau_{\kappa',T}^n$. We then have
\begin{eqnarray*}
&&\int_0^{\tau
_{\kappa',T}^n}e^{-\gamma s} |\check{\psi}^n(s)|\,ds\\
&&\qquad=\sum
_{l=1}^{r^n+1} \int_{\tau_{l-1}^n}^{\tau_l^n}e^{-\gamma s} |\check
{\psi}^n(s)|\,ds\\
&&\qquad=\sum_{l=1}^{r^n+1}\biggl(\int_{\tau
_{l-1}^n}^{(\tau_{l-1}^n+\eta^n)\wedge\tau_{l}^n}e^{-\gamma s}
|\check{\psi}^n(s)|\,ds+\int_{\tau_{l-1}^n+\eta^n}^{\tau_l^n\vee
(\tau_{l-1}^n+\eta^n) }e^{-\gamma s} |\check{\psi}^n(s)|\,ds
\biggr).
\end{eqnarray*}
The proof is now divided into three parts. We will show that,
under the conditions of the theorem,
%
%
\begin{eqnarray}\label{eqSSC1}
\Ex\Bigl[\sup_{1\leq l\leq
r^n+1}\sup_{\tau_{l-1}^n\leq s< (\tau_{l-1}^n+\eta^n)\wedge\tau
_l^n}|\check{\psi}^n(s)|\Bigr]&\leq& C\log^{k_0+2}n,
\\
\label{eqSSC2}
\Ex\Bigl[\sup_{1\leq l\leq r^n+1}\sup_{(\tau
_{l-1}^n+\eta
^n)\leq s< \tau_{l}\vee(\tau_{l-1}^n+\eta^n)}
|\check{\psi}^n(s)|\Bigr]&\leq& C\log^{k_0+2} n,
\end{eqnarray}
where we define $\sup_{(\tau_{l-1}^n+\eta^n)\leq s< \tau_{l}^n\vee
(\tau_{l-1}^n+\eta^n)}
|\check{\psi}^n(s)|=0$ if $\tau_l^n\leq\tau_{l-1}^n+\eta^n$.
Finally, we will show that
%
%
\begin{equation}\label{eqSSC3}
\Ex\biggl[\int_0^{\eta^n\wedge\tau_{\kappa
',T}^n}|\check
{\psi}^n(s)\,ds|\biggr]\leq C\log^{k_0}n.
\end{equation}
The proof of (\ref{eqSSC1}) hinges on the fact that, sufficiently
close to a change point~$\tau_l^n$, all the customer classes, $i$, for
which $h_i^n(\check{X}^n(s))=1$ for some $s$ in a~neighborhood of
$\tau_l^n$, will have similar values of $M_i^n(\check{X}^n(s))$. This
will follow from our gradient estimates for $\phi_{\kappa}^n$. The
proof of (\ref{eqSSC2}) hinges on the fact that,~$\eta^n$ time units
after a change point $\tau_l^n$ the queues of all the classes for
which $h^n(\check{X}^n(\tau_l^n))=0$ are small because, under the
tracking policy, these classes receive a significant share of the capacity.

Toward formalizing this intuition, define the following event on the
underlying probability space:
\begin{eqnarray*} \tilde{\Omega}(K)&=&\bigl\{|\check{X}^n(t)-\check
{X}^n(s)|\leq K \sqrt{n}\log^2 n(t-s)+K\log n,\\[-0.8pt]
&&\hspace*{135.3pt}
\mbox{ for all } s<t\leq\tau_{\kappa,T}^n\bigr\}\\[-0.8pt]
&&{}\cap_{\mH\subset\I} \biggl\{\mathcal{A}_{i,\mH
}^n(t)-\mathcal{A}_{i,\mH}^n(s)-\lambda_i^n(t-s)\geq\frac{\epsilon
_i }{2} n(t-s)-K\log n \\[-0.8pt]
&&\hspace*{183.5pt}\hspace*{12.1pt} \mbox{for all } s<t\leq\tau
_{\kappa,T}^n\biggr\}\\[-0.8pt]
&&{}\cap_{i\in\I} \{
|A_{i}^n(t)-A_{i}^n(s)-\lambda_i^n(t-s)|\leq\bar{\epsilon
}n(t-s)+K\log n \\[-0.8pt]
&&\hspace*{171.7pt}\mbox{for all } s<t\leq\tau
_{\kappa,T}^n\}.
\end{eqnarray*}
For each $0\leq t\leq\tau_{\kappa',T}^n$ and $i\in\I$ let
%
%
\begin{eqnarray}
\label{eqchecktaudefin}
\check{\varsigma}_i^n(t)&=&\sup\{s\leq t\dvtx h_i^n(\check{X}^n(s))=1\},
\\
\label{eqtildetaudefin}
\tilde{\varsigma}_i^n(t)&=&\inf\{
s\geq t\dvtx h_i^n(\check{X}^n(s))=1\}\wedge\tau_{\kappa',T}^n
\end{eqnarray}
and
%
%
\begin{equation}\label{eqydefin}
\hat{\varsigma}_i^n(t)=
\cases{\check{\varsigma}_i^n(t)+\eta^n,&\quad if
$Q_i^n(t)>4K\log n$,\cr
t, &\quad otherwise.}
\end{equation}
Then, we claim that on $\tilde{\Omega
}(K)$ and for all $t$ with $\tilde{\varsigma}_i^n(t)> \hat{\varsigma
}_i^n(t)$,
%
%
\begin{eqnarray}\label{eqSSCinterim}
&&\sup_{\hat{\varsigma}_i^n(t)\leq s< \tilde
{\varsigma}_i^n(t)}\bigl|\bigl(e\cdot\check{X}^n(s)\bigr)^+U_i^n(s)-\bigl(e\cdot\check
{X}^n(s)\bigr)^+h_i^n(\check{X}^n(s))\bigr|\nonumber\\[-8pt]\\[-8pt]
&&\qquad\leq12K\log n.\nonumber
\end{eqnarray}
Note that since $h_i^n(\cdot)\in\{0,1\}$, the
above is equivalently written as
%
%
\begin{equation}\label{eqSSCinterim2}
\sup_{\hat{\varsigma}_i^n(t)\leq s< \tilde
{\varsigma}_i^n(t)}Q_i^n(s)\leq12K\log n.
\end{equation}
In words,
when the process $\check{X}^n(t)$ enters a region in which
$h^n_i(\check{X}^n(\cdot))=0$ the queue of class $i$ will be drained
up to $12 K\log n$ within at most $\eta^n$ time units and it will
remain there up to $\tilde{\varsigma}_i^n(t)$. We postpone the proof
of (\ref{eqSSCinterim}) and use it in proceeding with the proof of
the theorem.

To that end, fix $l\geq0$ and let
\[
j^*_{l}=\min\mathop{\argmin}_{i\in\I}M_i^n(\check{X}^n(\tau_l^n)).
\]
Then, by the definition of the function $h^n$ in (\ref{eqixdefin}) we
have that $h_{j^*_l}^n(\check{X}^n(\tau_l^n))=1$ and $h_i(\check
{X}^n(\tau_l^n))=0$ for all $i\neq j^*_l$. In particular,
\begin{eqnarray*}
\check{\psi}^n(s)&=&\bigl(e\cdot\check
{X}^n(s)\bigr)^+h_{j^*_l}^n(\check{X}^n(s))M_{j^*_l}^n(\check{X}^n(s))
\\
&&{}-\sum_{i\in\I}Q_i^n(s)M_i^n(\check{X}^n(s))
\end{eqnarray*}
for all $s\in[\tau_{l}^n,(\tau_{l}^n+\eta^n)\wedge\tau_{l+1}^n)$.
Let
\[
\mathcal{J}(\tau_l^n)=\{i\in\I\dvtx Q_i^n(\tau_l^n-)>4K\log n\}.
\]
Then, simple manipulations yield
%
%
\begin{eqnarray}\label{eqcheckpsiinterim}
|\check{\psi}^n(s)|&\leq&\sum_{i\notin\mathcal{J}(\tau
_l^n)\cup\{j^*_l\}}Q_i^n(s)|M_i^n(\check{X}^n(s))|\nonumber\\
&&{}+
|M_{j_l^*}^n(\check{X}^n(s)|\biggl|\bigl(e\cdot\check
{X}^n(s)\bigr)^+-\sum_{i\in\mathcal{J}(\tau_l^n)\cup\{j^*_l\}
}Q_i^n(s)\biggr|\\
&&{}+\sum_{i\in\mathcal
{J}(\tau_l^n)}Q_i^n(s)|M_i^n(\check{X}^n(s))-M_{j_l^*}^n(\check
{X}^n(s))|.\nonumber
\end{eqnarray}
We turn to bound each of the elements on the right-hand side of (\ref
{eqcheckpsiinterim}). First, note that for all $i\notin\mathcal
{J}(\tau_l^n)\cup\{j^*_l\}$ it follows from (\ref{eqSSCinterim2}) that
\[
\sup_{\tau_{l}^n\leq s<(\tau_{l}^n+\eta^n)\wedge\tau
_{l+1}^n}Q_i^n(s)\leq12K\log n.
\]
Also, by (\ref{eqgradients1}) we have for all $i\in\I$ that
%
%
\begin{equation}\label{eqMboun4}
\sup_{0\leq s\leq\tau_{\kappa',T}^n}|M_i^n(\check{X}^n(s))|\leq
C\log^{k_1}n,\vadjust{\goodbreak}
\end{equation}
so that
%
%
\begin{equation}\label{eqitai1}
\sum_{i\notin\mathcal{J}(\tau_l^n)\cup\{j^*_l\}}Q_i^n(s)|M_i^n(\check
{X}^n(s))|\leq12IKC\log^{k_1+1}n
\end{equation}
for all $s\in
[\tau_{l}^n,(\tau_{l}^n+\eta^n)\wedge\tau_{l+1}^n)$ and a constant
$C$ that does not depend on $n$. From (\ref{eqSSCinterim2}) and from
the fact that $\sum_{i\in\I}Q_i^n(s)=(e\cdot\check{X}^n(s))^+$ we
similarly have that
%
%
\begin{equation}\label{eqitai2}
|M_{j_l^*}^n(\check{X}^n(s)|
\biggl|\bigl(e\cdot\check
{X}^n(s)\bigr)^+-\sum_{i\in\mathcal{J}(\tau_l^n)\cup\{j^*_l\}
}Q_i^n(s)\biggr|\leq12IC K\log^{k_1+1}n.\hspace*{-30pt}
\end{equation}
To
bound the last element on the right-hand side of (\ref
{eqcheckpsiinterim}) note that for each $i\in\mathcal{J}(\tau_l^n)$
there exists $\tau_l^n-\eta^n\leq t\leq\tau_{l}^n$ such that
$h_j^n(\check{X}^n(t))=1$. Otherwise, we would have a contradiction to
(\ref{eqSSCinterim}). We now claim that for each \mbox{$i\in\mathcal
{J}(\tau_l^n)$},
%
%
\begin{equation}\label{eqinterim747}
|M_i^n(\check{X}^n(s))-M_{j^*_l}^n(\check
{X}^n(s))|\leq\frac{C\log^{k_1+2} n}{\sqrt{n}}
\end{equation}
for all $s$ in $[\tau_l^n-\eta^n,\tau_l^n+\eta
^n]$. Indeed, by the definition of $\tilde{\Omega}(K)$, we have that
$|\check{X}^n(t)-\check{X}^n(s)|\leq C\log^{m+2}n$ for all $s,t$ in
$[\tau_l^n-\eta^n,\tau_l^n+\eta^n]$. As in the proof of (\ref
{eqgenbound}) [see, e.g., (\ref{eqMbound})] we have that
%
%
\begin{equation}\label{eqMbound2}
|M_i^n(x)-M_i^n(y)|\leq\frac{C\log^{k_2+m+2}n}{\sqrt{n}},\qquad i\in\I,
\end{equation}
for $x,y\in\mB_{\kappa'}^n$ with $|x-y|\leq
C\log^{m+2}n$. In turn,
%
%
\begin{equation}\label{eqMbound3}
|M_i^n(\check{X}^n(t))-M_i^n(\check
{X}^n(s))|\leq\frac{C\log^{k_2+m+2} n}{\sqrt{n}}=\frac{C\log
^{k_1+2} n}{\sqrt{n}}
\end{equation}
for all $i\in\I$ and all
$s,t\in[\tau_l^n-\eta^n,\tau_l^n+\eta^n]$ where we used the fact
that \mbox{$k_1=k_2+m$}. Since, for each $j\in\mathcal{J}(\tau_l^n)$, there
exists $\tau_l^n-\eta^n\leq t\leq\tau_{l}$ such that $h_j^n(\check
{X}^n(t))=1$ we have, by the definition of $h^n$ that $j\in\argmin
_{i\in\I} M_i^n(\check{X}^n(t))$ for such $t$ so that (\ref
{eqinterim747}) now follows from (\ref{eqMbound3}). Finally, recall
that\break $\sum_{i\in\I}Q_i^n(t)=(e\cdot\check{X}^n(s))^+\leq\kappa
\sqrt{n}\log^mn$ for all $s\leq\tau_{\kappa',T}^n$ and that
$k_0=k_1+m$ so that by~(\ref{eqinterim747})
\[
\sum_{i\in\mathcal{J}(\tau_l^n)}Q_i^n(s)|M_i^n(\check
{X}^n(s))-M_{j_l^*}^n(\check{X}^n(s))|\leq C\log^{k_0+2}n.
\]
Plugging this into (\ref{eqcheckpsiinterim}) together with (\ref
{eqitai1}) and (\ref{eqitai2}) we then have that, on $\tilde{\Omega}(K)$,
\[
\sup_{\tau_{l-1}^n\leq s<(\tau_{l-1}^n+\eta
^n)\wedge\tau_l^n}|\check{\psi}^n(s)| \leq C\log
^{k_1+m+2}n=CK\log^{k_0+2}n.
\]
This argument\vspace*{1pt} is repeated for each $l$. To complete the proof of (\ref
{eqSSC1}) note that, using (\ref{eqMboun4}) together\vspace*{-1pt} with ${\sup
_{0\leq s\leq\tau_{\kappa',T}^n}}|e\cdot\check{X}^n(s)|\leq\kappa
\sqrt{n}\log^mn$, we have that ${\sup_{0\leq s\leq\tau_{\kappa
',T}^n}}|\check{\psi}^n(s)|\leq C\sqrt{n}\log^{k_1+m}n$. Applying H\"
{o}lder's inequality we have that
\begin{eqnarray*}
&&\Ex\Bigl[\sup_{1\leq l\leq
r^n+1}\sup_{\tau_{l-1}^n\leq s< (\tau_{l-1}^n+\eta^n)}|\check{\psi
}^n(s)|\Bigr]\\
&&\qquad\leq
\Ex\Bigl[\sup_{1\leq l\leq r^n+1}\sup_{\tau_{l-1}^n\leq s< (\tau
_{l-1}^n+\eta^n)}|\check{\psi}^n(s)|1\{\tilde{\Omega}(K)\}
\Bigr]\\
&&\qquad\quad{}+
\Ex\Bigl[\max_{1\leq l\leq r^n+1}\sup_{\tau_{l-1}^n\leq s< (\tau
_{l-1}^n+\eta^n)}|\check{\psi}^n(s)|1\{\tilde{\Omega}(K)^c\}
\Bigr]\\
&&\qquad\leq
C\log^{k_0+2} n + C\sqrt{n}\log^{k_1+m}n C_1e^{-(C_2
{K}/{2})\log n}
\end{eqnarray*}
for redefined constants $C_1,C_2$ and (\ref{eqSSC1}) now follows by
choosing $K$ large enough.

We turn to prove (\ref{eqSSC2}). Rearranging terms in (\ref
{eqpsicheckdefin}) we write
\[
\check{\psi}^n(s) = \sum_{i\in\I}M_i^n(\check
{X}^n(s))\bigl(\bigl(e\cdot\check{X}^n(s)\bigr)^+h_i^n(\check{X}^n(s))-\bigl(e\cdot
\check{X}^n(s)\bigr)^+ U_i^n(s)\bigr),
\]
so that equation (\ref{eqSSC2}) now follows directly from (\ref
{eqSSCinterim}) and (\ref{eqMboun4}) through an application of H\"
{o}lder's inequality.

Finally, to establish (\ref{eqSSC3}), note that from the definition
of $\tau_{\kappa',T}^n$,
\begin{eqnarray*}
\sup_{0\leq t\leq\eta^n\wedge\tau_{\kappa',T}^n
} |\check{\psi}^n(s)|&\leq& I\sup_{0\leq t\leq\eta^n\wedge\tau
_{\kappa',T}^n}|\check{X}^n(t)|\sum_{i\in\I}M_i^n(\check{X}^n(t))
\\&\leq&
I\sup_{0\leq t\leq\eta^n\wedge\tau_{\kappa',T}^n}C\log
^{k_1}n|\check{X}^n(t)|\leq C \kappa\sqrt{n}\log^{k_1+m}n.
\end{eqnarray*}
In turn,
%
%
\begin{equation}\label{eqhowmanymoreinterims}
\Ex\biggl[\int_0^{\tau_0^n}e^{-\gamma t} |\check
{\psi
}^n(t)|\,dt\biggr]\leq C\log^{k_1+m}n=C\log^{k_0}n.
\end{equation}

We have thus proved (\ref{eqSSC1})--(\ref{eqSSC3}) and to conclude
the proof of the theorem it remains only to establish (\ref
{eqSSCinterim}). To that end, let $\check{\varsigma}_i^n(t)$,
$\tilde
{\varsigma}_i^n(t)$ and $\hat{\varsigma}_i^n(t)$ be as in~(\ref
{eqchecktaudefin})--(\ref{eqydefin}). Fix an interval $[l,s)\in
(\hat
{\varsigma}_i^n(t),\tilde{\varsigma}_i^n(t))$ such that
$Q_i^n(u)>2K\log n $ for all $u\in[l,s)$. By the\vadjust{\goodbreak} definition of the
tracking policy, (\ref{eqQtrack}) holds on this interval so that, on
$\omega\in\tilde{\Omega}(K)$,
%
%
\begin{eqnarray}\label{equntilwhen3}
Q_i^n(l)-Q_i^n(s)&\leq& \bar{\epsilon}n(t-s)
n-\frac{\epsilon_i}{2}n(t-s)+2K\log n\nonumber\\[-8pt]\\[-8pt]
&\leq& -\frac{\epsilon_i}{4}n(t-s)+2K\log n.\nonumber
\end{eqnarray}
Equation (\ref{eqSSCinterim}) now follows directly from (\ref
{equntilwhen3}). Indeed, note for all $t\leq\tau_{\kappa',T}$,
$Q_i^n(t)\leq(e\cdot\check{X}^n(t))^+\leq|\check{X}^n(t)|\leq
\kappa\sqrt{n}\log^mn$. Hence, $Q_i^n(\check{\varsigma}_i(t))\leq
\kappa\sqrt{n}\log^mn$. In turn, using (\ref{equntilwhen3}) and
assuming that $\tilde{\varsigma}_i^n(t)\geq\check{\varsigma
}_i^n(t)+\eta^n$ we have that $Q_i^n(\varsigma_{0,i}^n(t))\leq4K\log
n$ for some time $\varsigma_{0,i}^n(t)\leq\check{\varsigma
}_i^n(t)+\eta^n$ with $\eta^n$ as defined in (\ref{eqetandefin}).
Also, let
\[
\varsigma_{2,i}^n(t)=\inf\{t\geq\varsigma_{0,i}^n(t)\dvtx Q_i^n(t)\geq
12K\log n \}
\]
and
\[
\varsigma_{1,i}^n(t)=\sup\{t\leq\varsigma_{2,i}^n(t)\dvtx Q_i^n(t)\leq
8K\log n \}.
\]
Note that (\ref{equntilwhen3}) applies to any subinterval $[l,s)$ of
$[\varsigma_{1,i}^n(t),\varsigma_{2,i}^n(t))$. In turn, $\varsigma
_{2,i}^n(t)\leq\tilde{\varsigma}_i^n(t)$ would constitute\vspace*{1pt} a
contradiction to (\ref{equntilwhen3}) so that we must have that
$Q_i^n(s)\leq12K\log n$ for all $s\in[\varsigma_{0,i}^n(t),\tilde
{\varsigma}^n(t))$ with $\varsigma_{0,i}^n(t)\leq\tilde{\varsigma
}^n(t)+\eta^n$. Finally, note that $\varsigma_{0,i}^n(t)$ can be
taken to be $t$ if $Q_i^n(t)\leq4K\log n$.

This concludes the proof of (\ref{eqSSCinterim}) and, in turn, the
proof of the theorem.

\subsection*{\texorpdfstring{Proof of Lemma \protect\ref{lemafterstop}}{Proof of Lemma 6.2}} Let $T$,
$\tau_{\kappa',T}^n$ and $(x^n,q^n)$ be as in the statement of the
lemma. We first prove (\ref{eqafterstop2}). To that end, we claim
that, for all $T$ large enough,
%
%
\begin{equation}\label{eqthisisjustonemore}
\Ex_{x^n,q^n}\biggl[\int_{T\log
n}^{\infty} e^{-\gamma s} (e\cdot c)\bigl(e\cdot\check{X}^n(s)\bigr)^+\,ds
\biggr]\leq C\log^2n
\end{equation}
for some $C>0$ and all
$n\in\bbZ$. This is a direct consequence of Lemma 3 in~\cite{AMR02}
that, in our notation, guarantees that
\[
\Ex_{x^n,q^n}[|\check{X}^n(t)|]\leq C\bigl(1+|x^n|+ \sqrt{n}(t+t^2)\bigr)
\]
for all $t\geq0$ and some constant $C>0$. We use (\ref
{eqthisisjustonemore}) to prove Lemma \ref{lemafterstop}. The
assertion of the lemma will be established by showing that
\[
\Ex_{x^n,q^n}\biggl[\int_{\tau_{\kappa',T}^n}^{2T\log n} e^{-\gamma
s} (e\cdot c)\bigl(e\cdot\check{X}^n(s)\bigr)^+\,ds\biggr]\leq C\log^2n.
\]
To that end, applying H\"{o}lder's
inequality, we have
%
%
\begin{eqnarray}\label{eqinterim8}
&&\Ex_{x^n,q^n}\biggl[\int_{\tau_{\kappa
',T}^n}^{2T\log n} e^{-\gamma s} (e\cdot c)\bigl(e\cdot\check
{X}^n(s)\bigr)^+\,ds\biggr]\nonumber\\
&&\qquad\leq\Ex_{x^n,q^n}\Bigl[(2T\log
n-\tau_{\kappa',T}^n)^+\sup_{0\leq t\leq2T\log n}(e\cdot c)\bigl(e\cdot
\check{X}^n(t)\bigr)^+\,ds\Bigr]\nonumber\\[-8pt]\\[-8pt]
&&\qquad\leq
\sqrt{\Ex_{x^n,q^n}\bigl[\bigl((2T\log n-\tau_{\kappa
',T}^n)^+\bigr)^2\bigr]}\nonumber\\
&&\qquad\quad{}\times\sqrt{
\Ex_{x^n,q^n}\Bigl[\Bigl(\sup_{0\leq t\leq2T\log n}(e\cdot
c)\bigl(e\cdot\check{X}^n(t)\bigr)^+\,ds\Bigr)^2\Bigr]}.\nonumber
\end{eqnarray}

Using Lemma \ref{lemstrongappbounds} we have that
%
%
\begin{equation}\label{eqinterim9}
\Ex
_{x^n,q^n}\Bigl[\Bigl(\sup_{0\leq t\leq2T\log n}(e\cdot c)\bigl(e\cdot
\check{X}^n(t)\bigr)^+\,ds\Bigr)^2\Bigr]\leq Cn\log^6n
\end{equation}
for some $C>0$ (that can depend on $T$). Also, since
$m\geq3$,
\[
\Pd\{\tau_{\kappa',T}^n<2T\log n\}\leq\Pd
\Bigl\{\sup_{0\leq t\leq2T\log n}|\check{X}^n(t)|> \kappa'\sqrt{n}\log
^3 n -M\sqrt{n}\Bigr\}.
\]
Choosing $\kappa'$ (and in turn
$\kappa$ large enough) we then have, using Lemma \ref
{lemstrongappbounds}, that
%
%
\begin{equation}\label{eqtauprobbound}
\Pd\{\tau_{\kappa',T}^n<2T\log n\}
\leq\frac{C}{n^2}
\end{equation}
and hence, that
%
%
\begin{equation}\label{eqinterim10}
\Ex
_{x^n,q^n}\bigl[\bigl((2T\log n-\tau_{\kappa',T}^n)^+
\bigr)^2\bigr]\leq C.
\end{equation}
Plugging (\ref{eqinterim9}) and (\ref{eqinterim10}) into (\ref
{eqinterim8}) we
then have that
%
%
\begin{equation}\label{eqafterstopbound1}
\Ex_{x^n,q^n}\biggl[\int_{\tau_{\kappa
',T}^n}^{2T\log n}
e^{-\gamma s} (e\cdot c)\bigl(e\cdot\check{X}^n(s)\bigr)^+\,ds\biggr]\leq C\log
^2 n.
\end{equation}

To conclude the proof we will show that (\ref{eqafterstop1}) follows
from our analysis thus far. Indeed,
\begin{eqnarray*}
&&\Ex[e^{-\gamma\tau_{\kappa',T}^n}\phi
_{\kappa}^n(\check{X}^n(\tau_{\kappa',T}^n))]\\
&&\qquad\leq
\Ex_{x^n,q^n}^{U}\biggl[\int_{\tau_{\kappa',T}^n}^{2T\log
n}e^{-\gamma s} \sup_{0\leq s\leq2T\log n}(e\cdot c)\bigl(e\cdot\check
{X}^n(s)\bigr)^+ \,ds\biggr].
\end{eqnarray*}
The right-hand side here is bounded by $C\log^2n$ by the same argument
that leads to (\ref{eqafterstopbound1}).

\subsection*{\texorpdfstring{Proof of Lemma \protect\ref{lemmartingales}}{Proof of Lemma 6.3}} Recall that
$\check{W}^n$ is defined by $\check
{W}_i^n(t)=M_{i,1}^n(t)-M_{i,2}^n(t)$, where
\begin{eqnarray*}
M_{i,1}^n(t)&=&\mN_i^a(\lambda_i^n t)-\lambda_i^n t,
\\
M_{i,2}^n(t)&=&
\mN_i^d\biggl(\mu_i\int_0^t \bigl( \check{X}_i\lam(s)+\nu_i
n-U^n_i(s)\bigl( e\cdot\check{X}\lam(s)\bigr)^+ \bigr) \,ds \biggr)\\
&&{}-\mu
_i\int_0^t \bigl( \check{X}_i\lam(s)-U^n_i(s)\bigl( e\cdot\check
{X}\lam(s)\bigr)^+ \bigr) \,ds.
\end{eqnarray*}
The fact that each of the processes
$M_{i,1}^n(t)$ and $M_{i,2}^n(t)$ are square integrable martingales
with respect to the filtration $(\mathcal{F}_t^n)$ follows as in
Section 3 of \cite{pang2007martingale} and specifically as in Lemma
3.2 there.

Since, with probability 1, there are no simultaneous jumps of $\mathcal
{N}_i^a$ and $\mathcal{N}_i^d$, the quadratic variation process satisfies
\begin{eqnarray*}
[\check{W}_i^n]_t&=&[M_{i,1}^n]_t+[M_{i,2}^n]_t\\
&=&\sum_{s\leq
t}(\Delta M_{i,1}^n(s))^2+\sum_{s\leq t}(\Delta M_{i,2}^n(s))^2,
\end{eqnarray*}
where the last equality follows again from Lemma 3.1 in \cite
{pang2007martingale} (see also Example~5.65 in \cite{vandervaart}).
Finally, the predictable quadratic variation process satisfies
\begin{eqnarray*}
\langle\check{W}_i^n\rangle_t&=&\langle
M_{i,1}^n\rangle
_t+\langle M_{i,2}^n\rangle_t\\&=&\lambda_i^nt + \mu_i\int_0^t
\bigl( \check{X}_i\lam(s)+\nu_i n-U^n_i(s)\bigl( e\cdot\check{X}\lam
(s)\bigr)^+ \bigr) \,ds\\&=&
\int_0^t (\sigma_i^n(\check{X}^n(s),U^n(s)))^2\,ds,
\end{eqnarray*}
where the second equality follow again follows from Lemma 3.1 in \cite
{pang2007martingale} and the last equality from the definition of
$\sigma_i^n(\cdot,\cdot)$ [see (\ref{eqnbmdefn3})]. By Theorem
3.2 in~\cite{pang2007martingale} $((\check{W}_i^n(t))^2-[\check
{W}_i^n]_t,t\geq0])$ and $((\check{W}_i^n(t))^2-[\check
{W}_i^n]_t,t\geq0)$ are both martingales with respect to
$(\mathcal{F}_t^n)$. In turn, by the optional stopping theorem so are
the processes $\mathcal{M}_i^n(\cdot)$ and $\mathcal{V}_i^n(\cdot)$
as defined in the statement of the lemma. Finally, it is easy to verify
that these are square integrable martingales using the fact the time
changes are bounded for all finite $t$.
\end{appendix}


%
\printaddresses

\end{document}